\RequirePackage{fix-cm}
\documentclass[smallextended]{svjour3}       
\smartqed  
\usepackage{graphicx}
\usepackage[font=scriptsize]{caption}
\usepackage{amssymb,amsmath, cancel, lineno}
\usepackage{empheq}
\usepackage{url}
\usepackage{caption}[2008/04/01]
\usepackage{wrapfig}
\usepackage{subfig}
\usepackage{epstopdf}
\usepackage{nicefrac,lineno,subfloat}

\newcommand{\mc}{\mathcal}

\def\e1{{\varepsilon_{11}}}
\def\b1{{\beta_{11}}}
\def\bp3{{\beta_{33}}}
\def\ep3{{\varepsilon_{33}}}
\def\Re{{\rm Re \, }}
\def\mX{{\mathbb X}}
\def\Ltwo{{\mathbb L}^2 }

\usepackage{pifont}
\usepackage[tableposition=top]{caption}

\makeatletter
\def\elsartstyle{%
    \def\normalsize{\@setfontsize\normalsize\@xiipt{14.5}}
    \def\small{\@setfontsize\small\@xipt{13.6}}
    \let\footnotesize=\small
    \def\large{\@setfontsize\large\@xivpt{18}}
    \def\Large{\@setfontsize\Large\@xviipt{22}}
    \skip\@mpfootins = 18\p@ \@plus 2\p@
    \normalsize
}
\@ifundefined{square}{}{}
\makeatother

\MHInternalSyntaxOn
\providecommand*\phantomword[3][c]{%
\mathchoice
{\MT_phantom_word:NNnn #1\displaystyle {#2}{#3}}%
{\MT_phantom_word:NNnn #1\textstyle {#2}{#3}}%
{\MT_phantom_word:NNnn #1\scriptstyle {#2}{#3}}%
{\MT_phantom_word:NNnn #1\scriptscriptstyle {#2}{#3}}%
}
\def\MT_phantom_word:NNnn #1#2#3#4{%
\@begin@tempboxa\hbox{$\m@th#2#4$}%
\setlength\@tempdima{\widthof{$\m@th#2#3$}}%
\hbox{\hb@xt@\@tempdima{\csname bm@#1\endcsname}}%
\@end@tempboxa}
\MHInternalSyntaxOff

\newtheorem{rmk}{Remark}

\def\e{{\epsilon}}

\begin{document}

\title{Further stabilization  and exact observability results for voltage-actuated piezoelectric beams with magnetic effects
}


\author{Ahmet \" Ozkan \" Ozer      
}


\institute{
\at         Department of Mathematics \& Statistics \\
              University of Nevada-Reno, NV 89557\\
              Tel: 775-784-6774\\
              Fax: 775-784-6378\\
              \email{aozer@unr.edu}           
}

\date{Received: date / Accepted: date}

\maketitle

\begin{abstract}
It is well known that magnetic energy of the piezoelectric beam is relatively small,
and it does not change  the overall dynamics. Therefore,  the models, relying on electrostatic or quasi-static approaches, completely ignore the magnetic energy
stored/produced in the beam. A single piezoelectric beam model without the magnetic effects is known to be exactly observable and exponentially
stabilizable in the energy space. However, the model with the magnetic effects is proved to be not exactly observable / exponentially stabilizable in the energy space for almost all choices of material parameters. Moreover,  even strong stability is not achievable for many values of the material parameters. In this paper, it is shown that the uncontrolled system is exactly observable in a space larger than the energy space. Then, by using a $B^*-$type feedback controller,  explicit polynomial decay estimates are obtained for more regular initial data. Unlike the classical counterparts, this choice of feedback corresponds to the current flowing through the electrodes, and it matches better with the physics of the model. The results obtained in this manuscript have direct implications on the controllability/stabilizability of smart structures such as elastic beams/plates with piezoelectric patches and  the active constrained layer (ACL) damped beams/plates.

\end{abstract}
\section{Introduction}
\begin{figure}[h]
\centering
\vspace{-1.9\baselineskip}
\includegraphics[width=2.8in]{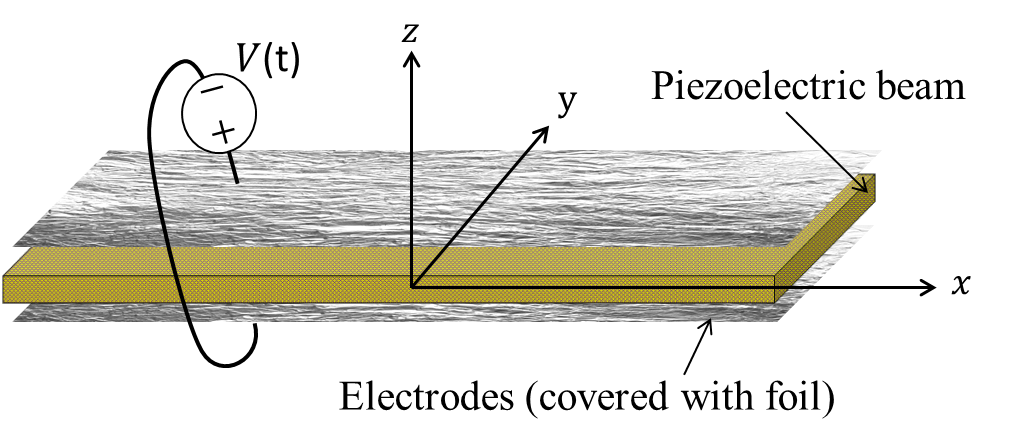}
\caption{ \small For a voltage-actuated beam/plate, when voltage $V(t)$ is supplied to the electrodes, an electric field is created between the electrodes, and therefore the beam/plate either shrinks or extends.}
\vspace{-1.9\baselineskip}
\label{pbeam}
\end{figure}

Piezoelectric material is an elastic beam/plate covered by electrodes at its top and bottom surfaces, insulated at the edges (to prevent fringing effects), and connected to an external electric circuit to create electric field between the top and the bottom electrodes (See Figure \ref{pbeam}).
It has a unique characteristic of converting mechanical energy to electrical and \emph{magnetic energy}, and vice versa. Therefore these materials could be used as both actuators or sensors. Moreover, since they are generally scalable, smaller, less expensive and more efficient than traditional actuators,  they have been employed in civil, industrial, automotive, aeronautic, and space structures.

In classical mechanics, it is very well known that equations of motion can be formulated either
through a set of differential equations, or through a variational principle, so-called
Hamilton's principle. In applying the Hamilton's principle, the functional is specified over a fixed time interval,
and the admissible variations of the generalized coordinates (independent variables) are taken to be
zero. The set of field equations for the piezoelectric beams/plates have been well established through
the coupling of beam/plate equations and Maxwell's equations. There are many different mathematical models proposed in the literature depending on the type of actuation; voltage, charge or current.

The linear models of piezoelectric beams incorporate three major effects and their interrelations: mechanical, electrical, and magnetic effects.  Mechanical effects are mostly modeled through Kirchhoff, Euler-Bernoulli, or Mindlin-Timoshenko small displacement assumptions. To include electrical and magnetic effects, there are mainly three approaches (due to Maxwell's equations): \emph{electrostatic, quasi-static,} and \emph{fully dynamic} \cite{Tiersten}.  Electrostatic approach is the most widely used  among the others. It completely excludes magnetic effects and their couplings with electrical and mechanical effects  (\cite{Banks-Smith,Dest,Rogacheva,Smith,Tiersten,Tzou,Yang}   and references therein). In this approach, even though  the mechanical equations are dynamic, electric field is not dynamically coupled. In other words, the electrical effects are assumed to be stationary. In the case of quasi-static approach
\cite{K-M-M,Yang}, magnetic effects are not completely ignored and electric charges have time dependence. The electromechanical coupling is still not dynamic though.

 Due to the small displacement assumptions, the stretching and bending motions of a single piezoelectric beam are completely decoupled. The bending equation without the electrical effects corresponds to the fourth order Euler-Bernoulli or Rayleigh/Kirchhoff beam equations; see i.e. \cite{O-M1}. Since the voltage control does not affect the bending equations at all,  we only consider the stretching equations in this paper. For a beam of length $L$ and thickness $h,$ the beam model (no damping) derived by Euler-Bernoulli small displacement assumptions, and electrostatic/quasi-static assumptions describe the stretching motion  as
\begin{subequations}
  \label{or}
\begin{empheq}[left={\phantomword[r]{0}{ }  \empheqlbrace}]{align}
  & \rho \ddot v-\alpha_1  v_{xx} = 0, &  (x,t)\in (0,L)\times \mathbb{R}^+ \\
  & v(0,t)= 0, ~ \alpha_1 v_x(L,t)= -\frac{ \gamma V(t)}{  h},&  t\in \mathbb{R}^+ \\
  &(v, \dot v)(x,0)=(v^0, v^1),&  x \in [0,L]
\end{empheq}
\end{subequations}
where  $\rho, \alpha_1, \gamma$ denote mass density,  elastic stiffness, and piezoelectric coefficients of the beam, respectively, $V(t)$ denotes the voltage applied at the electrodes, and $v$ denotes the longitudinal displacement of centerline of the beam. Throughout this paper, we use dots to denote differentiation with respect to time.

From the control theory point of view, it is well known that wave equation (\ref{or}) can be exactly controlled in the natural energy space (therefore the uncontrolled problem is exactly observable if the observability time is large enough). If we have the choice of a feedback in the form of a boundary damping $V(t)=- k~ \dot v(L,T)$ with $k>0,$  the solution of the closed-loop system is exponentially stable in the energy space (see, for instance \cite{Lions}).

 In the fully dynamic approach, magnetic effects are included, and therefore the wave behavior of the electromagnetic fields are accounted for, i.e. see \cite{O-M}.  These effects are experimentally observed to be minor on the overall dynamics for polarized ceramics  (see the review article \cite{Yang1}).  For a beam of length $L$ and thickness $h,$ the Euler-Bernoulli model with magnetic effects is derived in \cite{O-M1} as
\begin{subequations}
  \label{homo-vol}
\begin{empheq}[left={\phantomword[r]{0}{}  \empheqlbrace}]{align}
 \label{eq-main1} &\rho  \ddot v-\alpha   v_{xx}+\gamma\beta      p_{xx} = 0,   & \\
 \label{eq-main2}  & \mu  \ddot p   -\beta   p_{xx} + \gamma\beta  v_{xx}= 0, &  ~~ (x,t)\in (0,L)\times \mathbb{R}^+,
\end{empheq}
\begin{empheq}[left={\phantomword[r]{0}{ } \empheqlbrace}]{align}
\label{eq-main-ek0} & v(0)= p(0)=\alpha  v_{x}(L)-\gamma  \beta p_x(L)=0,&\\
\label{eq-main-ek}& \beta  p_x(L) -\gamma \beta v_x(L)= -\frac{V(t)}{h},&  t\in \mathbb{R}^+ \\
 \label{ivp}  &(v, p, \dot v, \dot p)(x,0)=(v^0,  p^0,  v^1,  p^1),& x \in [0,L]
 \end{empheq}
\end{subequations}
where  $\rho, \alpha, \gamma, \mu, \beta,$ and $V$ denote mass density per unit volume, elastic stiffness, piezoelectric coefficient, magnetic permeability, impermittivity coefficient of the beam, and voltage prescribed at the electrodes of the beam, respectively, and \begin{eqnarray}\label{alpha}\alpha=\alpha_1 + \gamma^2\beta.
 \end{eqnarray}
 Moreover,  $p=\int_0^x D_3(x,t) ~dt$ is the total charge at the point $x$ with $D_3(x,t)$ being the electric displacement along the transverse direction.
Observe that the term $\mu \ddot p$ in (\ref{homo-vol}) is due to the dynamic approach. If this term is ignored,   an elliptic-type differential equation is obtained, and once this equation is solved and back substituted to the mechanical equations, the system (\ref{homo-vol}) boils down to the system (\ref{or}) obtained in electrostatic and quasi-static approaches.

By using (\ref{alpha}), the boundary conditions (\ref{eq-main-ek}) can be simplified as the following
\begin{eqnarray} v_x(L)=-\frac{\gamma V(t)}{\alpha_1 h},\quad p_x(L)=-\frac{\alpha V(t)}{ \beta\alpha_1 h}.\label{eq-main3}\end{eqnarray}
The system (\ref{homo-vol}) with the simplified boundary conditions (\ref{eq-main3}) is a simultaneous controllability  problem with the control $V(t).$
Simultaneous controllability problems  were first introduced by \cite{Lions} and \cite{Russell}.
Controllability and stabilizability of the beam/plate with a control applied to \emph{a point/a curve} in the beam/plate cases   were investigated by a number of researchers including \cite{AvdoninII,AvdoninIII,Komornik-P2,Zuazua,Jaffard,Weiss-Tucsnak1}, and references therein.  By using a generalization of Ingham's inequality (with a weakened gap condition) (i.e., see \cite{Komornik-P})  and Diophantine's approximations \cite{Cassal}, exact controllability (observability)  in finite time, and stabilizability are obtained depending on the Diophantine approximation properties of the joints in the beam case, and how strategic the controlled curve is in the plate case.  Simultaneous controllability for general networks and trees are considered in \cite{Zuazua}. The controllability of a two interconnected beams (including the rotational inertia) by a point mass is considered in \cite{C-Z}. In this problem the weakened gap condition is a necessity. Notice that the system (\ref{homo-vol}) is a strongly coupled wave system whereas in \cite{Alabou1}-\cite{Alabou3} various other weakly coupled systems are considered. The methodology used in these papers is slightly different than ours. There are also research done in proving the controllability of various coupled parabolic systems, i.e. see \cite{T1,T2}, and \cite{T3}. The use of number theoretical results is unavoidable in \cite{T3}.

In this paper, we consider a coupled wave system  (\ref{homo-vol}) where the coupling terms  are at the order of the principal terms. The eigenvalues of the uncontrolled system ($V(t)\equiv 0$),  are all on the imaginary axis, and for almost all choices of parameters, they get arbitrarily close to each other (See Theorem \ref{close}). In other words, eigenvalues do not have a uniform gap.
Our first goal is to obtain the observability inequality for the uncontrolled system in a less regular space. Next, we choose a $B^*-$type feedback, i.e.   $V(t)=~\frac{1}{2h}\dot p(L)$ in (\ref{homo-vol}), to obtain the closed-loop system
\begin{subequations}
  \label{pdes-stabil}
\begin{empheq}[left={\phantomword[l]{}{ }  \empheqlbrace}]{align}
 \label{pdes-stabil-a} &\rho  \ddot v-\alpha   v_{xx}+\gamma  \beta    p_{xx} = 0   & \\
\label{pdes-stabil-b} & \mu \ddot p   -\beta   p_{xx} + \gamma  \beta  v_{xx}= 0, &  (x,t)\in (0,L)\times \mathbb{R}^+,
\end{empheq}
\begin{empheq}[left={\phantomword[l]{}{ }  \empheqlbrace}]{align}
 & v(0)= p(0)= \alpha  v_{x}(L)-\gamma \beta   p_x(L)=0, &\\
 & \beta  p_x(L) -\gamma \beta  v_x(L)= -\frac{ \dot p(L)}{2h^2},& t\in \mathbb{R}^+ \\
 &(v, p, \dot v, \dot p)(x,0)=(v^0, p^0, v^1, p^1),& x \in [0,L].
\end{empheq}
\end{subequations}
 In fact, the system (\ref{pdes-stabil}) is shown to be strongly stable \cite{cdcpaper}, but not exponentially stable  in the energy space for almost all choices of parameters  \cite{O-M1}. Based on the  observability inequality, we use the methods in \cite{Ammari-T} and \cite{Begout-Soria}  to obtain decay estimates for the solutions of the closed loop system (\ref{pdes-stabil}).  Notice that this type of feedback is very practical since it corresponds to the current flowing through the electrodes.

This paper is organized as follows. In section \ref{Sec-II}, we first prove that the uncontrolled system is well-posed in the interpolation spaces. In Section \ref{Sec-III},  we prove the exact observability results. In Section \ref{Sec-IV}, we give explicit decay rates for the solutions of the closed-loop system with the current feedback at the electrodes. Finally, in the Appendix, we briefly mention known results from number theory which are needed to prove our observability inequalities.

\section{Well-posedness}
\label{Sec-II}
The  energy associated with (\ref{homo-vol})  is given by
\begin{eqnarray}
\label{Energy-nat} \mathrm{E}(t)&=&\frac{1}{2}\int_0^L \left\{\rho  |\dot v|^2 +   \mu |\dot p|^2 + \alpha_1   |v_x|^2 + \beta\left| \gamma v_x- p_x\right|^2  \right\}~ dx,~~t\in \mathbb{R}.
\end{eqnarray}
We define the Hilbert space  \begin{eqnarray}
 \label{space3} && H^1_L(0,L)=\{v\in H^1(0,L): v(0)=0\}, \hspace{2em}\mX=(\Ltwo(0,L))^2
 \end{eqnarray}
and the complex linear space
\begin{eqnarray}\label{control}&& \mathrm{H} = \left(H^1_L(0,L)\right)^2 \times  \mX
 \end{eqnarray}
equipped with the energy inner product
{\small{\begin{eqnarray}
\nonumber && \left<\left( \begin{array}{l}
 u_1 \\
 u_2 \\
 u_3\\
 u_4
 \end{array} \right), \left( \begin{array}{l}
 v_1 \\
 v_2 \\
 v_3\\
 v_4
 \end{array} \right)\right>_{\mathrm{H}}= \left<\left( \begin{array}{l}
 u_3\\
 u_4
 \end{array} \right), \left( \begin{array}{l}
 v_3\\
 v_4
 \end{array} \right)\right>_{(\Ltwo(0,L)^2} +  \left<\left( \begin{array}{l}
 u_1 \\
 u_2
 \end{array} \right), \left( \begin{array}{l}
 v_1 \\
 v_2
 \end{array} \right)\right>_{\left(H^1_L(0,L)\right)^2}\\
\nonumber && :=\int_0^L \left\{\rho  u_3 \bar v_3 + \mu  u_4  \bar v_4\right\}~dx + \int_0^L \left\{\alpha_1   (u_1)_{x} (\bar v_1)_x +\beta\left( \gamma  (u_1)_x- (u_2)_x\right) \left( \gamma (\bar v_1)_x- (\bar v_2)_x\right) \right\}dx\\
\label{inner} &&= \int_0^L \left\{\rho  u_3 \bar v_3 + \mu  u_4  \bar v_4  + \left< \left( {\begin{array}{*{20}c}
   \alpha_1 + \gamma^2\beta  & -\gamma\beta\\
     -\gamma \beta & \beta  \\
\end{array}} \right) \left( \begin{array}{l}
 u_{1x} \\
 u_{2x}
 \end{array} \right), \left( \begin{array}{l}
  v_{1x} \\
  v_{2x}  \end{array} \right)\right>_{\mathbb{C}^2}\right\}~dx
 \end{eqnarray}}}
  where $\left<\cdot,\cdot\right>_{\mathbb{C}^2}$ is the inner product on $\mathbb{C}^2.$ Indeed, (\ref{inner}) is an inner product since the matrix $\left( {\begin{array}{*{20}c}
   \alpha_1 + \gamma^2\beta  & -\gamma\beta\\
     -\gamma \beta & \beta  \\
\end{array}} \right)$ is positive definite.

\subsection*{Interpolation spaces}
Define the operator $$A: {\text{Dom}}(A)\subset \mX\to \mX, \quad A= \left( {\begin{array}{*{20}c}
   -\frac{\alpha }{\rho}D_x^2  & \frac{\gamma\beta }{\rho}D_x^2  \\
       \frac{\gamma\beta}{\mu} D_x^2  &  -\frac{\beta}{\mu} D_x^2  \\
\end{array}} \right)$$
where
\begin{eqnarray}\label{dom-hom-A} {\rm {Dom}}(A)  = (H^2(0,L))^2~ \bigcap ~\{ (w_1, w_2)^{\rm T} \in (H^1_L(0,L))^2~:~ w_{1x}(L)= w_{2x}(L)= 0 \}.\quad~~
\end{eqnarray}
The operator $A$ can be easily shown to be a positive and self-adjoint operator, and since the $\rm{Dom} (A)$ is compactly embedded in $\mX,$ the operator $A^{-1}$ is compact, and therefore $A^{-1}$ has only countable many positive eigenvalues in its point spectrum, and all of its eigenvalues converge to zero. Therefore, the operator $A$ has has only countable many positive eigenvalues $\{\lambda_j\}_{j\in \mathbb N}$ in its point spectrum, and $|\lambda_j|\to \infty$ as $j\to \infty.$

Now we find the eigenvalues of A. Consider the eigenvalue problem
\begin{eqnarray}\label{e-A} A \left( \begin{array}{l}
 z_1\\
 z_2
 \end{array} \right)=\lambda \left( \begin{array}{l}
 z_1\\
 z_2
 \end{array} \right).\end{eqnarray}
Solving (\ref{e-A})  is equivalent to solving
\begin{subequations}
  \label{pdes-stabil10}
\begin{empheq}[left={\phantomword[r]{0}{ }  \empheqlbrace}]{align}
\label{st-010}  &\alpha   z_{1xx}-\gamma \beta  z_{2xx} = -\rho\lambda z_1   & \\
\label{st-020} &   \beta   z_{2xx} -\gamma \beta z_{1xx}= -\mu \lambda z_2 , &\\
\label{ivp-st-031}
 &z_1(0)=z_2(0)= z_{1x}(L)=z_{2x}(L)= 0.&
\end{empheq}
\end{subequations}
Define
\begin{eqnarray}
\label{lam1} \zeta_1 &=& \frac{1}{\sqrt{2}}\sqrt{\frac{\gamma ^2\mu }{\alpha_1}+\frac{\mu}{\beta}+\frac{\rho}{\alpha_1  }+\sqrt{\left(\frac{\gamma ^2\mu}{\alpha_1 }+\frac{\mu}{\beta}
+\frac{\rho}{\alpha_1  }\right)^2-\frac{4\rho  \mu}{\beta \alpha_1 }}} \\
\label{lam2} \zeta_2 &=& \frac{1}{\sqrt{2}}\sqrt{\frac{\gamma ^2\mu }{\alpha_1}+\frac{\mu}{\beta}+\frac{\rho}{\alpha_1  }-\sqrt{\left(\frac{\gamma ^2 \mu}{\alpha_1  }+\frac{\mu}{\beta}
+\frac{\rho}{\alpha_1  }\right)^2-\frac{4\rho  \mu}{\beta\alpha_1 }}}\\
\label{b-1} b_1 &=& \frac{1}{\gamma\mu}(\alpha_1\zeta_1^2-\rho)= \frac{1}{2}\left(\gamma+\frac{\alpha_1}{\gamma\beta}-\frac{\rho}{\gamma \mu}+\sqrt{\left(\gamma+\frac{\alpha_1}{\gamma\beta}-\frac{\rho}{\gamma \mu}\right)^2+\frac{4\rho}{  \mu}} \right)\\
\label{b-2} b_2 &=& \frac{1}{\gamma\mu}(\alpha_1\zeta_2^2-\rho)=\frac{1}{2}\left(\gamma+\frac{\alpha_1}{\gamma\beta}-\frac{\rho}{\gamma \mu}-\sqrt{\left(\gamma+\frac{\alpha_1}{\gamma\beta}-\frac{\rho}{\gamma \mu}\right)^2+\frac{4\rho}{  \mu}} \right).
\end{eqnarray}
Obviously, $\zeta_1,\zeta_2>0$ since
$$\left(\frac{\gamma ^2 \mu}{\alpha_1  }+\frac{\mu}{\beta}
+\frac{\rho}{\alpha_1  }\right)^2-\frac{4\rho  \mu}{\beta\alpha_1 }=\left(\frac{\gamma ^2 \mu}{\alpha_1  }+\frac{\mu}{\beta}
-\frac{\rho}{\alpha_1  }\right)^2+\frac{4\rho\mu\gamma^2}{\alpha_1^2 }>0,$$ and
 $$b_1, b_2\ne 0,\quad b_1\ne b_2, \quad  b_1 b_2=-\frac{\rho}{\mu}.$$

\begin{theorem} \label{main-thm-A} Let $ \sigma_j=\frac{(2j-1)\pi}{2L}, ~~j\in\mathbb{N}.$
The eigenvalue problem (\ref{e-A})  has distinct eigenvalues
 \begin{eqnarray}\label{EVs-A}\lambda_{1j}=\frac{\sigma_j^2}{\zeta_1^2},\quad \lambda_{2j}=\frac{\sigma_j^2}{\zeta_2^2},\quad  j\in \mathbb{N}\end{eqnarray}
 with the corresponding eigenfunctions
\begin{eqnarray}\label{evectors-A}  y_{1j}=\left( \begin{array}{c}
 1 \\
 b_1
 \end{array} \right)\sin \sigma_j x, && \quad y_{2j}=\left( \begin{array}{c}
 1 \\
b_2
 \end{array} \right)\sin \sigma_j x, \quad j\in \mathbb N.\end{eqnarray}
\end{theorem}
\textbf{Proof:} Using $\alpha=\alpha_1 + \gamma^2 \beta$  reduces (\ref{st-010}) and (\ref{st-020}) to
\begin{subequations}
  \label{pdes-stabil2}
\begin{empheq}[left={\phantomword[r]{0}{ }  \empheqlbrace}]{align}
\label{st-91}  &  z_{1xx}=\frac{-\lambda}{\alpha_1  }\left(\rho z_1 + \gamma \mu z_2\right)   & \\
\label{st-92} &  z_{2xx} =-\lambda\left(\frac{\gamma \rho}{\alpha_1  }z_1 +\left(\frac{\gamma ^2\mu }{ \alpha_1   }+\frac{\mu}{\beta} \right)  z_2\right). &
\end{empheq}
\end{subequations}
First, we find the eigenvalues of (\ref{EVs-A}). It is obvious that $\lambda = 0$ is not an eigenvalue since
the solution of (\ref{pdes-stabil2}) with (\ref{ivp-st-031}) is $z_1=z_2=0.$

We look for solutions of the form
\begin{eqnarray}z_{1j}=f_j \sin \sigma_{j} x, \quad z_{2j}=g_j \sin \sigma_{j} x.\label{modal}\end{eqnarray}
Solutions of this form satisfy all the homogeneous boundary conditions (\ref{ivp-st-031}). We seek
$f_j, g_j$ and $\lambda_j$  so that the system (\ref{pdes-stabil2}) is satisfied.
Substituting (\ref{modal}) into (\ref{pdes-stabil2}) we obtain
\begin{subequations}
  \label{pdes-stabil3}
\begin{empheq}[left={\phantomword[r]{0}{ }  \empheqlbrace}]{align}
\nonumber  &  \sigma_j^2 f_j =\frac{\lambda}{\alpha_1  }\left(\rho f_j + \gamma \mu g_j \right)   & \\
\nonumber  & \sigma_j^2 g_j =\lambda\left(\frac{\gamma \rho}{\alpha_1  } f_j+\left(\frac{\gamma ^2\mu }{ \alpha_1   }+\frac{\mu}{\beta} \right)g_j\right). &
\end{empheq}
\end{subequations}
The system above has nontrivial solutions if the following characteristic equation is satisfied
$$y_j^2 - \left(\frac{\gamma ^2 \mu }{\alpha_1 }+\frac{\mu}{\beta}+\frac{\rho}{\alpha_1  }\right)y_j + \frac{\rho \mu}{\beta \alpha_1  }=0.$$
where $y_j=\frac{\sigma_j^2}{\lambda}.$ Since $\left(\frac{\gamma ^2 \mu }{\alpha_1 }+\frac{\mu}{\beta}
+\frac{\rho}{\alpha_1  }\right)^2-\frac{4\rho  \mu }{\beta\alpha_1 }=\left(\frac{\gamma ^2\mu }{\alpha_1 }+\frac{\mu}{\beta}
-\frac{\rho}{\alpha_1  }\right)^2+\frac{4\rho \gamma ^2\mu }{\alpha_1  ^2}>0,$ a simple calculation shows that we have solutions
$ y_{j1}=\zeta_1^2, ~~y_{j2}=\zeta_2^2$ where $\zeta_1, \zeta_2\in\mathbb{R}$ are defined by (\ref{lam1}) and (\ref{lam2}), respectively.
Therefore $\lambda_{1j}=\frac{\sigma_j^2}{\zeta_1^2}, ~~~\lambda_{2j}=\frac{\sigma_j^2}{\zeta_2^2}, \quad j\in \mathbb{N},$ and (\ref{EVs-A}) follows.

Now we find the eigenvectors (\ref{evectors-A}). Let  $\lambda=\lambda_{1j}.$ Choosing $f_j=1$ yields  $g_j=b_1.$ The first eigenvector $y_{1j}$ follows from the solution $z_{1j}= \sin \sigma_j (x)$ and $z_{2j}=b_1\sin\sigma_j(x).$ Similarly, let $\lambda=\lambda_{2j}.$ Choosing $g_j=1$ yields $f_j=1/b_2.$ Hence the second eigenvector $y_{2j}$ follows from the solution $z_{1j}=\frac{1}{b_2} \sin\sigma_j(x)$ and $z_{2j}=\sin\sigma_j(x).$ $\square$

Obviously, the eigenvectors (\ref{evectors-A}) of $A$ are mutually orthogonal in $(H^1_L(0,L))^2$ by using the inner product defined by (\ref{inner}). Therefore, they form a Riesz basis in $(H^1_L(0,L))^2.$ Now we introduce the space $\mX_{\theta}={\rm Dom} (A^{\theta})$ for all $\theta\ge 0$ with the norm $\|\cdot \|_{\theta}=\|A^{\theta} \cdot\|_{\mX}.$ For example, using the definition of inner product $\left<\cdot, \cdot \right>_{(H^1_L(0,L))^2}$ in (\ref{inner}) yields
\begin{eqnarray}\nonumber \left<z_1, z_2\right>_{\mX_{1/2}}&=& \left<A^{1/2} z_1, A^{1/2} z_2\right>_{\mX}=\left<A z_1, z_2\right>_{\mX}=\left< z_1, z_2\right>_{(H^1_L(0,L))^2}.
\end{eqnarray}
The space $\mX_{-\theta}$ is defined to be the dual of $\mX_{\theta}$ pivoted with respect to $\mX.$ For example, the inner product on $\mX_{-1/2}$ is defined by
$$\left <z_1, z_2\right>_{\mX_{-1/2}}:=\left< A^{-1/2} z_1, A^{-1/2}z_2\right>_{\mX}=\left< A^{-1} z_1, z_2\right>_{\mX}.$$
Defining  $(H^1_L(0,L))^*$ to be the dual space of $H^1_L(0,L)$ pivoted with respect to $\Ltwo(0,L), $ we have
\begin{eqnarray}\label{duals}\mX_0=\mX, \quad \mX_{1/2}= (H^1_L(0,L))^2, \quad  \mX_{-1/2}=((H^1_L(0,L))^*)^2   \end{eqnarray}
Moreover, $\mX_1= {\rm Dom}(A)$ by the definition above. Note that the operator $A: \mX_{\theta}\to \mX_{\theta-1} $ can be boundedly extended or restricted for each $\theta\in\mathbb{R}.$

In fact, since the eigenvectors (\ref{evectors-A}) are mutually orthogonal  in $\mX_{\theta}$ for all $\theta\in \mathbb{R},$  every $U\in \mX_{\theta}$ has a unique expansion $U=\sum\limits_{k=1,2} \sum\limits_{j \in \mathbb{N}} c_{kj} y_{kj}$ where $c_{1j}, c_{2j}$ are complex numbers. Define the operator $A^{\theta}$ for all $\theta\in\mathbb{R}$ by
$$A^{\theta} U:=\sum\limits_{k=1,2} \sum\limits_{j \in \mathbb{N}} c_{kj} \lambda_{kj}^{\theta}y_{kj}.$$
Then
\begin{eqnarray}\|U\|_{\mX_{\theta/2}}^2=\left<A^{\theta} U, U\right>_{\mX}=\sum\limits_{k=1,2} \sum\limits_{j \in \mathbb{N}}\lambda_{kj}^{\theta}|c_{kj}|^2 \|y_{kj}\|^2_{\mX}.
\label{X_theta}
\end{eqnarray}
Similarly,
$$\|U\|_{\mX_{-\theta/2}}^2=\left<A^{-\theta} U, U\right>_{\mX}=\sum\limits_{k=1,2} \sum\limits_{j \in \mathbb{N}}\lambda_{kj}^{-\theta} |c_{kj}^2| \|y_{kj}\|^2_{\mX}.$$

\subsection*{Semigroup formulation}
Let $ \psi=(\psi_1, \psi_2, \psi_3, \psi_4)^{\rm T}=(v, p, \dot v, \dot p)^{\rm T}.$  Then  the system (\ref{homo-vol}) with the output $y(t)=\frac{1}{h}\dot p(L,t)$ can be put into the following state-space  formulation
\begin{subequations}
\label{Semigroup}
\begin{empheq}[left={\phantomword[l]{}{ }  \empheqlbrace}]{align}
&\dot \psi = \mc{A} \psi  + B V(t)= \left( {\begin{array}{*{20}c}
   0 & I_{2\times 2}  \\
   -A & 0  \\
\end{array}} \right) \psi + \left( \begin{array}{c}
0_{2\times 2}\\
B_0 \end{array} \right)V(t), &\\
&\psi(x,0) =  \psi ^0&\\
&y(t)=-B^*\psi = (0_{2\times 2}\quad B_0^*)~ \psi &
\end{empheq}
\end{subequations}
where \begin{eqnarray}
\nonumber  & B_0 \in \mathcal{L}(\mathbb{C}, \mX_{-1/2}), ~ \text{with} ~ B_0 V(t)= \left( \begin{array}{c}
0\\
-\frac{1}{h}\delta (x-L) \end{array} \right)V(t),& \\
\label{defb_0} & B_0^* \in \mc L( \mX_{1/2},  \mathbb{C}), ~ \text{with} ~  B^*\psi=(0_{2\times 2}\quad B_0^*)^{\rm T}\psi= -\frac{1}{h}\psi_4(L),&
\end{eqnarray}
 By the notation above we write $\mathrm H=\mX_{1/2}\times \mX.$ The operator $\mc A: {\rm {Dom}}(\mc A) \subset \mathrm{H}\to \mathrm{H}$ with the choice of the domain
\begin{eqnarray}\label{dom-hom} {\rm {Dom}}(\mc A)&=& \mX_1 \times \mX_{1/2}\\
  &=& (H^2(0,L))^2\times (H^1_L(0,L))^2 ~~\bigcap ~~\{\psi \in \mathrm{H}:~ \psi_{1x}(L)= \psi_{2x}(L)= 0 \}\quad\quad~~
\end{eqnarray}
is densely defined in $\mathrm{H}.$


 \begin{lemma} \cite{O-M1} \label{skew-adjoint}The infinitesimal generator $\mc{A}$  satisfies $\mc{A}^*=-\mc{A}$ on  $\mathrm{H},$ and  $\mc{A}$ and $\mc A^*$ are unitary, i.e.,
 \begin{eqnarray}\label{dang}\Re \left<\mc{A} \psi, \psi\right>_{\mathrm{H}}=\Re \left<\mc{A}^*\psi, \psi\right>_{\mathrm{H}}= 0.\end{eqnarray}
 Also, $\mc A$ has a compact resolvent.
\end{lemma}

Consider the uncontrolled system
\begin{subequations}
\label{Semigroup-H}
\begin{empheq}[left={\phantomword[l]{}{ }  \empheqlbrace}]{align}
&\dot \varphi = \mc{A} \varphi,&\\
& \varphi(x,0) =  \varphi ^0 &\\
&y(t)=-B^*\varphi. &
\end{empheq}
\end{subequations}

\begin{definition} The operator $B \in \mc L(\mathbb{C},\mathrm H_{-1})$ is an admissible control operator for $\{e^{\mc A t}\}_{t\ge 0}$ if
there exists a positive constant $c(T)$ such that for all $u \in H^1 (0,T )$,
$$\left\|  \int_0^T  e^{\mc A (T- t)} B u (t) dt \right\|_{\mathrm H} \le c(T) \| u \|_{\Ltwo (0,T)}. $$
\end{definition}

\begin{definition} The operator $B^*\in \mc L({\rm Dom}(\mc A), \mathbb{C})$ is an admissible observation operator for $\{e^{\mc A^*t} \}_{t\ge 0}$ if
there exists a positive constant $c(T)$ such that for all $\varphi^0\in {\rm Dom}(\mc A)$
$$\int_0^T \|B^* e^{\mc A^* t} \varphi^0\|^2~ dt \le c(T) \|\varphi^0\|^2_{\mathrm H}.$$
\end{definition}

The operator $B^*$ is an admissible observation operator for $\{e^{\mc A^*t} \}_{t\ge 0},$ if and only if $B$ is an admissible control operator for $\{e^{\mc A t}\}_{t\ge 0}$ \cite[pg. 127]{Weiss-Tucsnak}).

It is proved in \cite{O-M1} that both $B$ and $B^*$ operators are admissible. Now the theorem  on well-posedness of (\ref{Semigroup}) is now immediate.

\begin{theorem} \cite{O-M1} \label{w-p}
Let $T>0,$ and $V(t)\in \Ltwo(0,T).$ For any $\psi^0 \in \mathrm{H},$ there exists  positive constants $c_1(T),c_2(T)$ and a unique solution  to (\ref{Semigroup}) with $\psi \in C ([0,T]; \mathrm{H} ),$ and
      \begin{eqnarray}\label{conc}\|\psi\|^2_{\mathrm{H}} &\le& c_1 (T)\left\{\|\psi^0\|^2_{\mathrm{H}} + \|V(t)\|^2_{\Ltwo(0,T)}\right\},\\
      \label{conc-a} \|y(t)\|^2_{\Ltwo(0,T)} ~dt &\le& c_2(T) \left\{\|y(0)\|_{\mathrm H}^2 + \|V(t)\|^2_{\Ltwo(0,T)}\right\}.
      \end{eqnarray}
\end{theorem}

 We have the following theorem characterizing the eigenvalues and eigenfunctions of $\mc A.$

\begin{theorem} \label{main-thm} Let $\sigma_j=\frac{(2j-1)\pi}{2L},~~  j\in \mathbb{N}.$   The eigenvalue problem $\mc AY=\lambda Y$ has distinct eigenvalues
 \begin{eqnarray}\label{EVs} &\tilde \lambda_{1j}^\mp= \mp i \sqrt{\lambda_{1j}}=\frac{\mp i\sigma_j}{\zeta_1},\quad \tilde \lambda_{2j}^\mp=\mp i \sqrt{\lambda_{2j}}=\frac{\mp i\sigma_j}{\zeta_2},\quad  j\in \mathbb{N}&
 \end{eqnarray}
Since $\tilde \lambda_{1j}^-=-\tilde \lambda_{1j}^+, \quad \tilde \lambda_{2j}^-=-\tilde \lambda_{2j}^+, \quad j\in\mathbb{N},$ the corresponding eigenfunctions are
\begin{eqnarray}\nonumber Y_{1j}=\left( \begin{array}{c}
 \frac{1}{\tilde\lambda_{1j}^+} \\
 \frac{b_1}{\tilde\lambda_{1j}^+} \\
 1\\
b_1
 \end{array} \right)\sin \sigma_j x, &&~ Y_{-1j}=\left( \begin{array}{c}
 \frac{1}{\tilde\lambda_{1j}^+} \\
\frac{ b_1}{\tilde\lambda_{1j}^+} \\
-1 \\
-b_1
 \end{array} \right)\sin \sigma_j x,\\
\label{evectors}  ~Y_{2j}=\left( \begin{array}{c}
 \frac{1}{\tilde\lambda_{2j}^+} \\
\frac{ b_2}{\tilde\lambda_{2j}^+} \\
1\\
b_2
 \end{array} \right)\sin \sigma_j x, &&~Y_{-2j}=\left( \begin{array}{c}
 \frac{1}{\tilde\lambda_{2j}^+} \\
 \frac{b_2}{\tilde\lambda_{2j}^+} \\
-1\\
-b_2
 \end{array} \right)\sin \sigma_j x,\quad j\in \mathbb N\end{eqnarray}
 where  $\zeta_1, \zeta_2, b_1$ and $b_2$ are defined by (\ref{lam1})-(\ref{b-2}).
  The function
\begin{eqnarray}
\label{sol} \varphi(x,t)&=& \sum\limits_{j \in \mathbb{N}}  \left[c_{1j}Y_{1j}e^{\tilde\lambda_{1j}^+}t
 + d_{1j} Y_{-1j} e^{-\tilde\lambda_{1j}^+t} + c_{2j}Y_{2j} e^{\tilde\lambda_{2j}^+t}+d_{2j}Y_{-2j}e^{-\tilde\lambda_{2j}^+t}\right]\quad\quad\quad
\end{eqnarray}
solves  (\ref{Semigroup-H}) for the initial data
 \begin{eqnarray}\nonumber \varphi^0&=&\sum\limits_{j \in \mathbb{N}}  \left[c_{1j}Y_{1j}
 + d_{1j} Y_{-1j}  + c_{2j}Y_{2j} +d_{2j}Y_{-2j}\right]\\
 \label{init-F} &=&\sum\limits_{j \in \mathbb{N}} \left( \begin{array}{c}
\frac{1}{\tilde\lambda_{1j}^+}(c_{1j}+ d_{1j})+ \frac{1}{\tilde\lambda_{2j}^+}(c_{2j} + d_{2j}) \\
 \frac{b_1}{\tilde \lambda_{1j}^+} (c_{1j}+ d_{1j})+\frac{b_2}{\tilde\lambda_{2j}^+} (c_{2j} + d_{2j}) \\
(c_{1j}- d_{1j})+(c_{2j} -d_{2j})\\
b_1 (c_{1j}- d_{1j})+b_2(c_{2j} - d_{2j})
 \end{array} \right)\sin\sigma_j x
\end{eqnarray}
where $\{c_{kj}, d_{kj}, \quad k=1,2, \quad j\in\mathbb{N}\}$ are complex numbers such that
\begin{eqnarray}\label{norm-eq}&\|\varphi^0\|_{\mathrm{H}}^2 \asymp \sum\limits_{j \in \mathbb{N}} \left(|c_{1j}|^2 + |d_{1j}|^2+ |c_{2j}|^2 + |d_{2j}|^2\right), ~~{\rm i.e.,}&\\
\label{norm-init} & \tilde C_1~ \|\varphi^0\|_{\mathrm{H}}^2 \le \sum\limits_{j \in \mathbb{N}} \left(|c_{1j}|^2 + |d_{1j}|^2+ |c_{2j}|^2 + |d_{2j}|^2\right)\le \tilde C_2 ~\|\varphi^0\|_{\mathrm{H}}^2 &\end{eqnarray}
with two positive constants $\tilde C_1, \tilde C_2$ which are independent of the particular choice of $\Psi^0\in \mathrm H.$

\end{theorem}

\textbf{Proof:} Let $W=(W_1, W_2)^{\mathrm T}.$ Solving the eigenvalue problem $\mc A W= \tilde \lambda W$ is equivalent to solving $A W_1= -\tilde \lambda^2 W_1$ and $W_2=\tilde \lambda W_1.$ Since $\{\lambda_{1j}, ~\lambda_{2j}, j\in\mathbb{N}\}$ defined by (\ref{EVs-A}) are the eigenvalues of $A,$ it follows that $\tilde \lambda_{1j}^\mp=\mp i \sqrt{\lambda_{1j}}$ and $\tilde \lambda_{2j}^\mp=\mp i \sqrt{\lambda_{2j}},~~j\in\mathbb{N},$ and therefore (\ref{EVs}) follows. (\ref{sol}) and (\ref{init-F}) follow from (\ref{EVs}),(\ref{evectors}) and Theorem \ref{main-thm-A}. For the proof of  (\ref{norm-init}), see \cite{O-M1}. $\square$

It is easy to show that the eigenfunctions $\{Y_{kj},\quad k=-2,-1,1,2, \quad j\in \mathbb{N}\}$ are mutually orthogonal in $\mathrm H$ (with respect to the inner product (\ref{inner})). Therefore, they form a Riesz basis in $\mathrm H.$ This result also follows from the fact that we have a skew-symmetric operator $\mc A$ with a compact resolvent (see Lemma \ref{skew-adjoint}).

For $\theta\in \mathbb{R},$ we define   the space
\begin{equation}\label{stars-A}\mc S_{\theta}:=\{\sum\limits_{k =1,2}\sum\limits_{j\in \mathbb{N}} c_{kj}Y_{kj} + d_{kj}Y_{-kj}~:~ \sum\limits_{k =1,2}\sum\limits_{j\in \mathbb{N}}{|\tilde \lambda_{kj}|^{2\theta}}\left(\left|c_{kj}\right|^2 + \left|d_{kj}\right|^2\right)< \infty\}\end{equation}
by the completion of eigenvectors $\{Y_{kj},\quad k=-2,-1,1,2, \quad j\in \mathbb{N}\}$ with respect to the norm
 \begin{eqnarray} \left\|U\right\|^2_{\mc S_{\theta}}&=& \sum\limits_{k =1,2}\sum\limits_{j\in \mathbb{N}}{|\tilde \lambda_{kj}|^{2\theta}} \left( \left|c_{kj}\right|^2 + \left|d_{kj}\right|^2\right).\label{4stars-A}\end{eqnarray}

\begin{rmk} \label{equiv} For the simplicity of the calculations in the next sections, we use  the equivalent norm $\left\|U\right\|_{\mc S_{\theta}}=\left(\sum\limits_{k =1,2}\sum\limits_{j\in \mathbb{N}}{| 2j-1 |^{2\theta}} \left( \left|c_{kj}\right|^2 + \left|d_{kj}\right|^2\right)\right)^{\frac{1}{2}}.$ This follows from $\zeta_1,\zeta_2>0.$
\end{rmk}

Denote the space $\mc S_{-\theta}$ to be dual of $\mc S_{\theta}$ pivoted with respect to $\mc S_0:=\mathrm H=(H^1_L(0,L))^2\times (\Ltwo(0,L))^2.$ By (\ref{X_theta})
 \begin{eqnarray} \nonumber \mc S_{1}&=&\mX_1 \times \mX_{1/2}={\rm Dom}(\mc A)\\
   \nonumber \mc S_{0}&=& \mX_{1/2} \times \mX =\mathrm H\\
\nonumber   \nonumber \mc S_{-1}&=& \mX \times \mX_{-1/2}\end{eqnarray}
 Let $0<\varepsilon<\frac{1}{2}.$ By (\ref{stars-A}), we can also define the interpolation spaces
 $$ \mX_{1/2+\varepsilon/2}=[\mX_1, \mX_{1/2}]_{1-\varepsilon/2}, \mX_{\varepsilon}=[\mX_{1/2}, \mX]_{1-\varepsilon/2}$$ so that
 $$[\mc S_1, S_0]_{1-\varepsilon}=\mc S_{\varepsilon}=  \mc \mX_{1/2+\varepsilon/2} \times \mX_{\varepsilon/2},$$   $$[\mc S_2, S_1]_{1-\varepsilon}=\mc S_{1+\varepsilon}=\mX_{1+\varepsilon/2} \times \mX_{1/2+\varepsilon/2}$$ and their duals $ \mc S_{-1-\varepsilon}$ and $\mc S_{-\varepsilon}$ pivoted with respect to $\mc S_0=\mathrm H;$ see \cite{Triebel} for more information on interpolation spaces. We have the following dense compact embeddings
\begin{eqnarray}\nonumber \mc S_{1+\varepsilon} \subset \mc S_1\subset \mc S_{\varepsilon} \subset \mc S_0 \subset \mc S_{-\varepsilon} \subset \mc S_{-1}\subset \mc S_{-1-\varepsilon}.\end{eqnarray}
 With the notation above $\mc S_{-1-\varepsilon}=\mX_{-\varepsilon/2}\times  \mX_{-1/2-\varepsilon/2}.$

Now we have the following result from  \cite{Weiss-Tucsnak}:

Since $A:\mX_{1}\to \mX$ can be uniquely extended (or restricted) to $\tilde {\tilde A}: \mX_{\theta}\to \mX_{\theta-1}$ for any $\theta\in\mathbb{R},$ the infinitesimal generator $\mc{A}: \mathrm S_1 \to \mc S_0$ can be uniquely extended to $\tilde{\tilde{\mc{A}}}: \mc S_{-\varepsilon} \to  \mc S_{-1-\varepsilon}.$

 \begin{corollary}\label{ext2} The semigroup $\{e^{\mc A t}\}_{t\ge 0}$ with the generator $\mc{A}: \mathrm S_1 \to \mc S_0$ has a unique extension to a  contraction semigroup    $\{e^{\tilde{\tilde{\mc A }}t}\}_{t\ge 0}$    on $\mc S_{-1-\varepsilon}$  with the generator $\tilde{\tilde{\mc{A}}}: \mc S_{-\varepsilon} \to  \mc S_{-1-\varepsilon}$ for any $0<\varepsilon<\frac{1}{2}.$
\end{corollary}

\section{Exact observability}
\label{Sec-III}

We start with the definition of exact observability.

\begin{definition} The pair  $(\mc A^*, B^*)$ is exactly observable in time $T>0$  if there exists a positive constant $C(T)$ such that for all $\varphi^0\in\mathrm H$
$$\int_0^T \|B^* e^{\mc A^* t} \varphi^0\|^2~ dt \ge C(T) \|\varphi^0\|^2_{\mathrm H}.$$
\end{definition}
The following  theorem is proved in \cite{O-M1}.



 \begin{theorem} \label{close} Assume that $\frac{\zeta_2}{\zeta_1}\in \mathbb{R}-\mathbb{Q}.$ The eigenvalues $\{\tilde \lambda_{1j}^\mp=\frac{\mp i\sigma_j}{\zeta_1},\quad \tilde \lambda_{2m}^\mp=\frac{\mp i \sigma_m}{\zeta_2}, ~~j,m\in \mathbb{N}\}$ given by Theorem \ref{main-thm}  can get arbitrarily close to each other  for some choices of $j$ and $m.$ Therefore, the system (\ref{Semigroup-H}) is not exactly observable on $\mathrm{H}$.
 \end{theorem}

For the system (\ref{Semigroup-H}), Ingham-type theorems (see i.e. \cite{Komornik-P}, \cite{Weiss-Tucsnak}) can not be used to obtain the observability inequality since they require a uniform gap between the eigenvalues. This type of problem is well studied for joint structures with a point mass at the joint ( see \cite{Komornik-P} and references therein), or  for networks of strings/beams with different lengths (see \cite{Zuazua} and references therein). The main idea of proving observability result is based on the use of divided differences \cite{Ulrich}, the generalized Beurling's theorem, and the Diophantine's approximation.  We try the idea in \cite{Komornik-P} with the following technical result to prove our main observability result.


\begin{lemma} \label{irr3} Assume that $\frac{\zeta_2}{\zeta_1}\in \mathbb{R}-\tilde{\mathbb{Q}}$ where the set $\tilde{\mathbb{Q}}$ is defined in Theorem \ref{irr2}. Then there exists a number $\tilde\tau>0$ such that if
\begin{eqnarray}0<|\tilde\lambda_{kj}^+-\tilde\lambda_{lm}^+|\le\tilde\tau,\quad k,l=1,2,\quad j,m\in\mathbb{N}\label{eq20}
\end{eqnarray} then $k\ne l$ and
 $$|\tilde\lambda_{1j}^+-\tilde\lambda_{2m}^+|\ge \frac{C_\alpha}{|\tilde \lambda_{1j}^+|^{\alpha}},\quad |\tilde\lambda_{1j}^+-\tilde\lambda_{2m}^+|\ge \frac{C_{\alpha}}{|\tilde \lambda_{2m}^+|^{\alpha}} $$
 for every $\alpha> 1,$ with a constant $C_\alpha$ independent of the particular choice of $\tilde\lambda_{1j}$ and $\tilde\lambda_{2m}.$
\end{lemma}

\textbf{Proof of Lemma \ref{irr3}:} Since $\zeta_1,\zeta_2\in \mathbb{R}-\tilde {\mathbb{Q}},$ we have $\tilde\lambda_{kj}^+\ne \tilde\lambda_{lm}^+$ for any $k,l=1,2, ~~j,m\in \mathbb{N}.$ If we choose $\tilde\tau< \left(\frac{\pi }{L}\right) \mathop {\min } \left(\frac{1}{\zeta_1}, \frac{1}{\zeta_2}\right),$ (\ref{eq20}) is satisfied. This implies that $k\ne l.$ By Theorem \ref{irr2}, there exists a sequence of odd integers $\{\tilde  p_j\}, \{\tilde q_j\}\to \infty$ and $\alpha>1$ such that
$$\left|\tilde q_j\frac{\zeta_2}{\zeta_1}- \tilde p_j \right|\ge \frac{\tilde C_\alpha}{(\tilde q_j)^\alpha}.$$
Therefore
$$\left|\tilde\lambda_{1j}^+-\tilde\lambda_{2m}^+\right|=\frac{\pi}{2L}\left| {\frac{(2j+1)}{\zeta_1}}-{\frac{(2m+1) }{\zeta_2}}\right|\ge \frac{\pi}{2L}\frac{\tilde C_\alpha}{(2j+1)^\alpha}\ge \frac{\tilde C_\alpha}{|\tilde \lambda_{1j}|^{\alpha}},$$
and there is always a rational number $r$ such that $(2j+1)=r(2m+1)$ so that $C_\alpha$ can be chosen smaller to get
$$\left|\tilde\lambda_{1j}^+-\tilde \lambda_{2m}^+\right|\ge \frac{C_\alpha}{|\tilde \lambda_{2m}^+|^\alpha}. \square$$

We also need the following technical lemma from  \cite[Chap. 9]{Komornik-P} which is a slightly different version of the result obtained in \cite{Ulrich}:

\begin{lemma}\label{ana}Given an increasing sequence $\{s_{n}\}$ of real numbers satisfying
\begin{equation}s_{n+2}-s_{n}\ge 2\tau\quad {\rm for ~ all ~} n,\label{gap}\end{equation}
fix  $0<\tau'\le \tau$ arbitrarily and introduce the divided differences of $\{e_n(t), e_{n+1}(t)\}$ of exponential functions $\{e^{is_nt},e^{is_{n+1}t}\}$ by
\begin{eqnarray}\label{divided}e_n(t)=e^{is_nt},\quad e_{n+1}(t)=\frac{e^{is_{n+1}t}-e^{is_{n}t}}{s_{n+1}-s_n}.\end{eqnarray}
 Then  there exists positive constants $\tilde c_3(T)$ and $\tilde c_4(T)$ such that  $T>\frac{2\pi}{\tau}$
 $$\tilde c_3(T)  \sum\limits_{n =  - \infty }^\infty  |a_n|^2\le \int_0^T |f(t)|^2 \ dt \le \tilde c_4(T) \sum\limits_{n =  - \infty }^\infty  |a_n|^2$$ holds for all functions  given by the sum
$f(t)=\sum\limits_{n =  - \infty }^\infty  a_n e_n(t)~:~ \sum\limits_{n =  - \infty }^\infty{\left| {a_n } \right|^2 }<\infty. $
\end{lemma}

Now we can prove our main observability result:

\begin{theorem}\label{main} Let $\frac{\zeta_2}{\zeta_1}\in \mathbb{R}-\tilde {\mathbb{Q}}$  and  $T > 2L (\zeta_1+\zeta_2).$ Then there exists a constant $C=C(T)>0$ such that solutions $\varphi$ of the problem  (\ref{Semigroup-H}) satisfy the following observability estimate:
\begin{eqnarray}
\label{obs} & \int_0^T | B^*\varphi|^2 ~dt \ge  C(T) \|\varphi^0\|_{\mc S_{-1-\varepsilon}}^2.&
\end{eqnarray}
where $\mc S_{-1-\varepsilon}$ is defined by (\ref{4stars-A}).
\end{theorem}

\textbf{Proof:}  Let $s_{1j}=\frac{\sigma_j}{\zeta_1}=\frac{(2j-1)\pi}{2L\zeta_1}$ and $s_{2j}=\frac{\sigma_j}{\zeta_2}=\frac{(2j-1)\pi}{2L\zeta_2}$ for $j\in\mathbb{N}.$ The set of eigenvalues (\ref{EVs}) can be rewritten as
\begin{eqnarray}\label{forgap}  \tilde \lambda_{kj}^{\mp}=\mp is_{kj}, ~~k=1,2, ~~ j\in \mathbb{N}.
\end{eqnarray}

Since $\mc A^*=-\mc A,$ the function $\varphi=e^{\mc A^*t}\varphi^0,$ given explicitly by (\ref{sol}), solves (\ref{Semigroup-H}), and by (\ref{defb_0}) and (\ref{EVs})-(\ref{norm-init})
$$B^*\varphi=\sum\limits_{k =1,2}\sum\limits_{j\in \mathbb{N}} b_k (-1)^j\left(c_{kj} e^{is_{kj}^+t} +  d_{kj}e^{-is_{kj}^+t}\right). $$
By (\ref{4stars-A}), showing (\ref{obs}) is equivalent to showing
\begin{eqnarray}
\nonumber \int_0^T |B^*\varphi|^2 ~dt &=& \int_0^T \left|\sum\limits_{k =1,2}\sum\limits_{j\in \mathbb{N}} b_k (-1)^j\left(c_{kj} e^{is_{kj}^+t} +  d_{kj}e^{-is_{kj}^+t}\right)  \right|^2~dt\\
\label{suck21} &\ge & C(T) \sum\limits_{k =1,2}\sum\limits_{j\in \mathbb{N}}\frac{\left|c_{kj}\right|^2 + \left|d_{kj}\right|^2}{|\tilde\lambda_{kj}|^{2+2\varepsilon}}.
\end{eqnarray}
Let's rearrange $\{\mp s_{kj}^+:~ k=1,2, ~ j\in \mathbb{N}\}$ into an increasing sequence of $\{s_n, ~n\in \mathbb{N}\}.$ Denote the coefficients $\{(-1)^jb_k c_{kj}, (-1)^j b_k d_{kj}\}$ by $g_n$ (recall that $b_1,b_2\in \mathbb{R}-\{0\}$).  Then showing (\ref{suck21}) is equivalent to showing
\begin{eqnarray}
\label{suck22} \int_0^T |B^* \varphi|^2 ~dt &=& \int_0^T \left|\sum\limits_{n\in \mathbb{N}} g_n e^{is_n t}\right|^2~dt\ge  C(T) \sum\limits_{n\in \mathbb{N}} \frac{|g_n|^2}{|s_{n}|^{2+2\varepsilon}}.
\end{eqnarray}

 Let $n^+(r)$ denotes the largest number of terms of the sequence $\{s_n, ~n\in \mathbb{N}\}$ contained in an interval of length $r.$ Then
$$\frac{L (\zeta_1+\zeta_2) ~r }{\pi}-1\le n^+(r)\le \frac{L (\zeta_1+\zeta_2) ~r }{\pi}+1$$
Therefore  $D^+=\mathop {\lim }\limits_{r \to \infty } \frac{n^+(r)}{r}=\frac{L (\zeta_1+\zeta_2) }{\pi}.$ Now let $\tau=\frac{\pi}{2L}{\rm min} \left(\frac{1}{\zeta_1}, \frac{1}{\zeta_2}\right)$  so that  \begin{equation}\label{suck1}s_{n+2}- s_{n}\ge 2\tau, \quad {\rm for~ all} ~~n.\end{equation}
Note that the condition $T>\frac{2\pi}{\tau}$ can be replaced by $T>2\pi D^+=2L(\zeta_1+\zeta_2) $ (See Prop. 9.3 in \cite{Komornik-P}). Now we fix $0<\tau'<\tau$ and define sets $A_1$ and $A_2$  of integers by
\begin{eqnarray}
\nonumber A_1&:=&\{m, m+1 \in \mathbb{N}~:~s_{m+1}-s_{m}< \tau'\}\\
\nonumber A_2&:=&\{k \in \mathbb{N}~:~|s_{k}- s_{m}|\ge \tau', ~~ m\in A_1\}.
\end{eqnarray}
Observe that index $n$ of  the eigenvalues $\{s_n\}$ belongs to either $A_1$ or $A_2.$
For $m\in A_1,$ the exponents $\{s_m,  s_{m+1}\}$ form a \emph{chain of close exponents} for $\tau'$ and there is no chain of close exponents longer than two elements. For $m\in A_1,$ the divided differences $e_m(t), e_{m+1}(t)$ of the exponential functions are defined by (\ref{divided}).
Therefore, by Lemma \ref{ana} for all $T>2L(\zeta_1+\zeta_2) $ we have
\begin{equation}\int_0^T \left|\sum\limits_{n \in\mathbb{N}} a_n e_{n}(t)\right|^2~dt\asymp \sum\limits_{n \in\mathbb{N}}{\left| {a_n } \right|^2 }.\label{suck11}\end{equation}
 If $m\in A_1,$ we rewrite the sums as
 $$ \sum\limits_{n =  m }^{m+1}{g_ne^{i s_n t} }= \sum\limits_{n =  m }^{m+1}{a_n e_n(t)}$$
 where
  $a_{m}=g_m+\frac{a_{m+1}}{s_{m+1} - s_m},\quad a_{m+1}= g_{m+1}(s_{m+1}-s_m).$
  Then there exists a constant $C>0$ independent of $m$ such that
\begin{equation}\label{suck10} \sum\limits_{n =  m }^{m+1}{|g_n|^2 |s_{m+1}-s_{m}| }^2\le  C \sum\limits_{n =  m }^{m+1}{|a_n|^2}.\end{equation}
 By  Lemma \ref{irr3}, there exists a constant $C_\alpha>0$ such that
 $$|s_{m+1}-s_m|^2\ge \frac{C_{\alpha}}{|s_m|^{2\alpha}}, ~~{\rm and}~~|s_{m+1}-s_m|^2\ge \frac{C_{\alpha}}{|s_{m+1}|^{2\alpha}} $$
 where $\alpha> 1.$ Therefore by (\ref{suck10}) for all $\alpha= 1+\varepsilon$
 $$\sum\limits_{n =  m }^{m+1}{\frac{|g_n|^2}{|s_n|^{2+2\varepsilon}} } \le \frac{C}{C_{1+\varepsilon}} \sum\limits_{n =  m }^{m+1}{|a_n|^{2}}. $$
 On the other hand, if $n\in A_2,$ with the choice of a smaller $C_{1+\varepsilon}$ (if necessary) we get
 \begin{equation}\label{suck100}{\frac{|g_n|^2}{ |s_{n}|^{2+2\varepsilon} }}\le \frac{C}{C_{1+\varepsilon}}  {|g_n|^2},\end{equation}
  and (\ref{suck11}),(\ref{suck10}), and (\ref{suck100})  imply that for $T>2\pi D^+=2L(\zeta_1+\zeta_2)$
  $$\sum\limits_{n \in\mathbb{N} }{\frac{ |g_n|^2}{|s_n|^{2+2\varepsilon}}} \le \frac{C}{C_{1+\varepsilon}}\sum\limits_{n \in\mathbb{N} }{|a_n|^{2}}.$$
  This together with (\ref{suck11}) implies (\ref{suck22}), and therefore (\ref{obs}) holds. $\square$

\begin{corollary}\label{cor} Let $\frac{\zeta_2}{\zeta_1}\in {\tilde{\tilde {\mathbb{Q}}}}$  and  $T > 2L (\zeta_1+\zeta_2).$ Then there exists a constant $C=C(T)>0$ such that solutions of the problem  (\ref{Semigroup-H}) satisfy the following observability estimate:
\begin{eqnarray}
\label{obs10} & \int_0^T |B^*\varphi|^2 ~dt \ge  C(T) \|\varphi^0\|_{\mc S_{-1}}^2.&
\end{eqnarray}
where $\mc S_{-1}$ is defined by (\ref{stars-A}).
 \end{corollary}

 \textbf{Proof:} If we replace the inequality of  (\ref{aptal})  by (\ref{yeni}), then the proofs of Lemma \ref{irr3} and  Theorem \ref{main} can be adapted for $\varepsilon=0.$ This implies that the observability inequality (\ref{obs}) holds as $\mc S_{-1+\varepsilon}$ is replaced by $\mc S_{-1}.$ $\square$

\begin{rmk} Note that the lower bound of the control time $T=2L (\zeta_1+\zeta_2)$ obtained in Theorem \ref{main} and Corollary \ref{cor} is optimal.
The optimality of the control time can be  obtained by using the theory (i.e. see \cite{AvdoninII,AvdoninIII}). However, since the main scope of the paper is proving the polynomial stability and investigating the decay rates, we plan to use their idea in the upcoming research of exact controllability of the elastic beam/patch system. 
\end{rmk}
\section{Stabilization}
\label{Sec-IV}
The signal $\dot p(L, t)$ is the observation dual to the control
operator $B$ in (\ref{Semigroup}), and so we choose the feedback $V (t) = -\frac{1}{2} B^* z=  \frac{1}{2h} \dot p(L, t)$ in (\ref{Semigroup}). Also, since $\dot p(L, t)$ is the total
current at the electrodes, this variable can be measured easier than the velocity of the beam at one end.
The system (\ref{pdes-stabil}) can be put in the following form
\begin{subequations}
\label{Sg}
\begin{empheq}[left={\phantomword[r]{0}{ }  \empheqlbrace}]{align}
&\dot z(t) = \mc A_d  z(t) = \left( {\begin{array}{*{20}c}
   0 & I_{2\times 2}  \\
   -A & -\frac{1}{2}B_0 B_0^*  \\
\end{array}} \right)z, &\\
&z(x,0) =  z^0,&\\
\label{out}&y(t)=-B^*z(t)&
\end{empheq}
\end{subequations}
where $ \mc A_d : \text{Dom} ( \mc A_d ) \subset \mathrm H \to \mathrm H$ and $\text{Dom} ( \mc A_d )$ is defined by
 \begin{eqnarray}\nonumber &{\rm {Dom}}(\mc A_d)=\left\{z \in (H^2(0,L))^2\times (H^1_L(0,L))^2~:~ z_1(0)=z_2(0)=0,\right.&~~\\
\label{sem-dom} &\left.\alpha  z_{1x}(L)-\gamma \beta z_{2x}(L)=0, ~~\beta  z_{2x}(L) -\gamma \beta z_{1x}(L)= -\frac{ z_4(L)}{2h^2} \right\}.&
\end{eqnarray}
Note that the system above is equivalent to the system studied in \cite{Ammari-T}.

\begin{definition} The semigroup $\{e^{\mc A_d t}\}_{t\ge 0}$ with the generator $\mc A_d$ is exponentially stable on $\mathrm H$ if there exists constants $M, \mu>0$ such that $\|e^{\mc A_d t}\|_{\mathrm H}\le M e^{-\mu t}$ for all $t\ge 0.$
\end{definition}

\begin{theorem} \cite{cdcpaper,O-M1} \label{stronglystable}
\begin{itemize}
  \item [(i)] $ \mc A_d  : {\rm {Dom}}(  \mc A_d ) \to \mathrm{H}$ is the infinitesimal generator of a $C_0-$semigroup of contractions.
Therefore for every $T\ge 0,$ and $z^0 \in \mathrm H,$ $z$ solves (\ref{Sg}) with
$ z\in  C\left([0,T]; \mathrm{H}\right).$
\item [(ii)] The spectrum $\sigma( \mc A_d )$ of $ \mc A_d $ has all isolated eigenvalues. The semigroup $\{e^{ \mc A_d t}\}_{t\ge 0}$ is strongly stable on $\mathrm{H}$ if and only if  $\frac{\zeta_1}{\zeta_2}\ne \frac{2n+1}{2m+1},$ for some $ n,m \in \mathbb{N}.$  where $\zeta_1$ and $\zeta_2$ are defined by (\ref{lam1}) and (\ref{lam2}), respectively.
      \item [(iii)] Assume that $\frac{\zeta_2}{\zeta_1}\in \mathbb{R}-\mathbb{Q}.$  The semigroup $\{e^{ \mc A_d t}\}_{t\ge 0}$ is not exponentially stable on $\mathrm{H}.$
\end{itemize}
\end{theorem}

\subsection*{\textbf{Decay estimates}}

We need the following results to prove our main stabilization results given by Theorem \ref{stab1}.

\begin{lemma}  \label{Ammari's} \cite[Lemma 4.4]{Ammari-T}  Let $\{\mc E_k\}_{k\in \mathbb{N}}$ be a sequence of real numbers satisfying
$$\mc E_{k+1}\le \mc E_k - C \mc E_{k+1}^{2+\alpha}$$
where $C>0$ and $\alpha>-1$ are constants. Then there exists a positive constant $M=M(\alpha, C)$ such that
$$\mc E_k\le \frac{M}{(k+1)^{\frac{1}{1+\alpha}}}, \quad k\in \mathbb{N}.$$
\end{lemma}

\begin{lemma} \cite[Theorem 2.2]{Begout-Soria} \label{Begout} Let $m\in\mathbb{R}^+,$  $\omega_1 : (m,\infty) \to (0,\omega_1(m))$ and  $\omega_2 : (m, \infty) \to (0,\infty)$  be convex and increasing, and convex and decreasing functions, respectively, with $\omega_1(\infty)=0$ and $\omega_2(\infty)=\infty.$ Let $\Phi: (0,\omega_1(m))\to (0,\infty)$ and $\Psi: (0,\omega_2(m))\to (0,\infty)$ be concave and increasing functions with $\Phi(0)=0,$  $\Psi(\infty)=\infty,$ and for all $t>m$
$$1\le \Phi(\omega_1(t)) \Psi(\omega_2(t)).$$
 Then for  $j\in\mathbb{N},$ $j\ge m,$ and any $0 \ne f=\{f_j\}_{j\in \mathbb{N}}\in l_1 (\mathbb{N};\mathbb{R}),$ we have
$$1\le \Phi\left(\frac{ \sum\limits_{j\in \mathbb{N}} |f_j | \omega_1(j)}{\sum\limits_{j\in \mathbb{N}} |f_j|}\right)\Psi\left(\frac { \sum\limits_{j\in \mathbb{N}} |f_j| \omega_2(j)}{\sum\limits_{j\in \mathbb{N}} |f_j|}\right)$$
where $\{f_j \omega_1{j}\}_{j\in \mathbb{N}}, \{f_j \omega_2 (j)\}_{j\in \mathbb{N}}\in l^1(\mathbb{N};\mathbb{R})$, and therefore
\begin{eqnarray}\sum\limits_{j\in \mathbb{N}} |f_j| \omega_1(j) \ge \mc G_{\Phi,\Psi}^{-1} \left(\frac {\sum\limits_{j\in \mathbb{N}} |f_j| \omega_2(j)}{\sum\limits_{j\in \mathbb{N}} |f_j|}\right), \quad \mc G_{\Phi,\Psi}(j)=\frac{1}{\Psi^{-1}\left(\frac{1}{\Phi(j)}\right)}.\label{salak}
\end{eqnarray}
\end{lemma}

 Lemma \ref{Begout} is the discrete version of the H\"{o}lder's inequality originally proved in \cite{Begout-Soria}. That is, we use the discrete measure $\mu$  with the measurable weights $\omega_1$ and $\omega_2.$ For instance $\int_{m}^\infty |f(x)| d\mu (x)=\sum\limits_{(m \le) j\in \mathbb{N}} |f(x_j) | \omega_1(j)$ where $f=\{f_j\}_{(m\le)j\in\mathbb{N}}\in l^1(\mathbb{N};\mathbb{R}).$

Now we are ready to prove our main stabilization result:

\begin{theorem} \label{stab1}
  \noindent (I) Suppose that $\frac{\zeta_2}{\zeta_1}\in \tilde {\mathbb Q}$ where $\tilde {\mathbb Q}$  is defined in Theorem \ref{irr2}. Then for all $t\ge 0,$ there exists a positive constant $M_1$ such that
      \begin{eqnarray}\label{sonuc0}\|z(t)\|_{\mathrm H}^2\le \frac{M_1}{(t+1)^{\frac{1}{1+\varepsilon}}}\|z^0\|^2_{\mc S_{1+\varepsilon}}.\end{eqnarray}
  \noindent (II) Suppose that $\frac{\zeta_2}{\zeta_1}\in \tilde{\tilde {\mathbb Q}}$ where $\tilde{\tilde {\mathbb Q}}$  is defined in (\ref{yeni}). Then for all $t\ge 0,$ there exists a positive constant $M_2$ such that
      \begin{eqnarray}\|z(t)\|_{\mathrm H}^2\le \frac{M_2}{t+1}\|z^0\|^2_{\mc S_{1}}.\label{sonuc1}\end{eqnarray}
\end{theorem}

\textbf{Proof:} 
 Assume that $\psi$ and $\varphi$ solve (\ref{Semigroup}) and (\ref{Semigroup-H}) with the initial data $\psi^0= 0,$ $\varphi^0=z^0,$ and with $V(t)=B^* z$ so that $z=\varphi + \psi$ solves (\ref{Sg}). By  (\ref{obs}) we have
 \begin{eqnarray}
\label{obs-fake} & \int_0^T | B^* \varphi|^2 ~dt \ge  C(T) \| z^0\|_{\mc S_{-1-\varepsilon}}^2.&
\end{eqnarray}
   On the other hand since $B^* z= B^* \varphi+ B^* \psi,$ we can write
   $$|B^* \varphi|\le |B^*  z| + |B^* \psi|.$$
By (\ref{conc})  with $V(t)=B^*z,$ $\psi^0=0,$ and $y(t)=B^*\psi$
$$|B^* \psi|\le  |B^*  z|,$$
and by (\ref{obs-fake}) we obtain
   \begin{eqnarray}
\label{obs-fake2} & \int_0^T | B^* z|^2 ~dt \ge  C(T) \|z^0\|_{\mc S_{-1-\varepsilon}}^2.&
\end{eqnarray}
This proves the observability result for (\ref{Sg}).

To apply Lemma \ref{Begout}, we choose
  $$m=1/2,\quad \omega_1(j)=\frac{1}{(2j-1)^{2+2\varepsilon}},\quad \omega_2(j)=(2j-1)^2,$$
  and two functions $\Phi(t)$ and $\Psi(t):$
  $$\Phi(t)=\frac{1}{\omega_1^{-1}(t)}=\frac{2}{\left(\frac{1}{t}\right)^{\frac{1}{2+2\varepsilon}}+1}, \quad \Psi(t)=\omega_2^{-1}=\frac{\sqrt{t}+1}{2}.$$ Then $\mc G_{\Phi,\Psi}(t)=\frac{1}{t^{\frac{1}{1+\varepsilon}}}$ and $\mc G_{\Phi,\Psi}^{-1}(t)=\frac{1}{t^{1+\varepsilon}}.$  Denoting $\{|f_j|\}$ and $\{\tilde \lambda_{kj}\}$ by $\{|c_{kj}|^2+|d_{kj}|^2\} $ and $\{\tilde \lambda_{j}\},$ respectively,
  \begin{equation}\sum\limits_{j\in \mathbb{N}} |f_j | ~\omega_1(j)=\|f\|_{\mc S_{-1-\varepsilon}}^2, \quad \sum\limits_{j\in \mathbb{N}} |f_j | =\|f\|_{\mathrm H}^2, \quad \sum\limits_{j\in \mathbb{N}} |f_j | ~\omega_2(j)=\|f\|_{\mc S_{1}}^2\label{angut}\end{equation}
  where we used (\ref{stars-A}) and Remark \ref{equiv}.   By using (\ref{salak}) together with (\ref{angut})   we obtain

\begin{eqnarray}\|z^0\|_{\mc S_{-1-\varepsilon}}^2\ge \|z^0\|_{\mathrm H}^2~ \mc G_{\Phi,\Psi}^{-1}  \left(\frac {\|z^0\|_{\mc S_1}^2}{\|z^0\|^2_{\mathrm H}}\right)=\|z^0\|^2_{\mathrm H}\left(\frac{\|z^0\|^2_{\mathrm H}}{\|z^0\|_{\mc S_1}^2}\right)^{1+\varepsilon}=\frac{\|z^0\|^{4+2\varepsilon}_{\mathrm H}}{\|z^0\|_{\mc S_1}^{2+2\varepsilon}}.\label{angut2}\end{eqnarray}
By Theorem \ref{stab1}, (\ref{obs-fake2}),  (\ref{angut2}), and the fact that the function $t \mapsto  \|z(t)\|_{\mathrm H}$ is nonincreasing, we obtain
\begin{eqnarray}\nonumber \|z(T)\|_{\mathrm H}^2&=&\|z^0\|_{\mathrm H}^2- \int_0^T |B^* z|^2~dt \\
\nonumber &\le& \|z^0\|_{\mathrm H}^2- C(T)\|z^0\|_{\mc S_{-1-\varepsilon}}^2\\
\nonumber &\le& \|z^0\|_{\mathrm H}^2-C(T)\frac{\|z^0\|^{4+2\varepsilon}_{\mathrm H}}{\|z^0\|_{\mc S_1}^{2+2\varepsilon}}\\
\label{salak2} &\le& \|z^0\|_{\mathrm H}^2-C(T)\frac{\|z(T)\|^{4+2\varepsilon}_{\mathrm H}}{\|z^0\|_{\mc S_1}^{2+2\varepsilon}}.\end{eqnarray}
The estimate (\ref{salak2}) remains valid in successive intervals $[mT, (m+1)T].$ So, for all $m\ge 0,$ we have
\begin{eqnarray}\label{salak3}\|z((m+1) T)\|_{\mathrm H}^2\le  \|z(m T)\|_{\mathrm H}^2-C( T)\frac{\|z((m+1) T)\|^{4+2\varepsilon}_{\mathrm H}}{\|z^0\|_{\mc S_1}^{2+2\varepsilon}}.\end{eqnarray}
By letting $\mc E_m=\frac{\|z(m T)\|^{2}_{\mathrm H}}{\|z^0\|_{\mc S_1}^{2}}, $ (\ref{salak3}) gives
$$\mc E_{m+1}\le \mc E_m -C( T) (\mc E_{m+1})^{2+\varepsilon}, \quad m\in \mathbb{N}.$$
Hence, by Lemma \ref{Ammari's}, there exists a constant $M_1=M_1(C(T))>0$ such that (\ref{sonuc1}) follows.

 The proof of (II) is similar to the proof of (I) modulo a few simple changes. We take $\varepsilon=0,$ and use the observability inequality (\ref{obs10}) instead of (\ref{obs}). $\square$


\section{Conclusion and Future Research}
 The main result of this paper is to show that magnetic effects
in piezoelectric beams, even though small, have a dramatic effect on exact observability
and stabilizability. The piezoelectric beam model, without magnetic effects, is exactly
observable and exponentially stabilizable, by  a $B^*-$type feedback  However, when magnetic effects
are included, the beam is not exactly observable or exponentially stabilizable. By the $B^*-$ type feedback, the beam can be exactly observable and polynomially stabilizable for the initial data $z^0$ in $\mathcal S_1$ and $\mc S_{1+\varepsilon}$ when the ratio $\frac{\zeta_2}{\zeta_1}$ is in the sets $ \mathrm{\tilde Q}$ or $\tilde{\tilde {\mathrm Q}},$ respectively. These sets are of Lebesque measure zero even though they are uncountable.

\setlength\columnsep{\marginparsep}
\begin{wrapfigure}{r}[-.1cm]{0.45\textwidth}
\label{sand}
\vspace{-1.9\baselineskip}
\subfloat[An elastic beam with piezoelectric patches where the voltages $V_T(t)$ and $V_B(t)$ are applied to the top and bottom piezoelectric patches, repsectively.]{\label{lbl5}
  \includegraphics[width=0.42\textwidth]{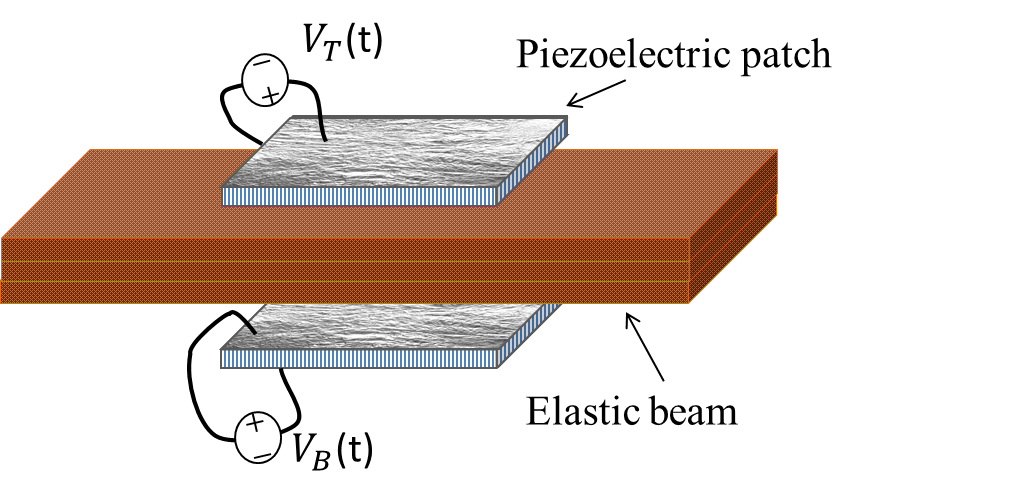}}\par
 \subfloat[Active constrained layer (ACL) damped beam/plate where the voltage $V_T$ is applied to the piezoelectric patch.]{\label{lbl3}
 \includegraphics[width=0.42\textwidth]{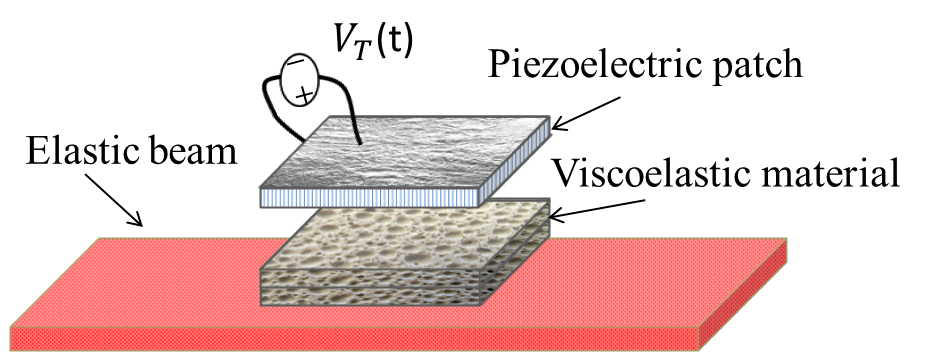}}
\vspace{-0.9\baselineskip}
\end{wrapfigure}

A single piezoelectric beam model using the Euler Bernoulli or Mindlin-Timoshenko small displacement assumptions is assumed to contract/extend only (by the linear theory). The voltage control does not even affect the bending motions \cite{accpaper}. A related and more interesting problem is to find the decay rates of the elastic beam/patch system (see Figure \ref{lbl5}). Once the magnetic effects are included \cite{accpaper}, the behavior of the system differs substantially from the classical counterparts \cite{Banks-Smith,J-T,Tucsnak-a,Tucsnak-b} which use electrostatic or quasi-static assumptions. In this model, the stretching equations (\ref{homo-vol}) are coupled to the bending (and rotation) equations, and it is similar in nature to the transmission problem proposed by Lions \cite{Lions}. The beam domain is divided into three sub-domains; first and the third for the pure elastic and the second for piezo-elastic coupling.  Previous research on controllability of elastic beam/plate with piezoelectric patches  without
magnetic effects showed that the location of the patch(es) on the beam/plate
strongly determines the controllability and stabilizability. This paper, \cite{cdcpaper,O-M1},   and \cite{accpaper}
suggest that the controllability and stabilizability depends on not only the location of the
patches but also the system parameters. This is currently under investigation.

 Our results in this paper also have strong implications on the controllability of smart sandwich structures such as Active Constrained Layer (ACL) damped  structures (see Figure \ref{lbl3}). The classical sandwich beam or plate is an engineering model for a  multi-layer beam consisting of  ``face" plates and ``core" layers that are orders of magnitude are more compliant than the face plates. ACL damped beams are sandwich structures of elastic, viscoelastic,  and piezoelectric layers. These structures are being successfully used for a variety of applications such as spacecraft, aircraft, train and car structures, wind turbine blades, boat/ship superstructures. i.e. see \cite{Baz}. The modeling and control strategies developed in \cite{Hansen3,O-Hansen1,O-M,O-M1,O-Hansen3,O-Hansen4}  play a key role in accurate analysis of these structures.  The controllability/stabilization problems in the case of  voltage actuation is still an open problem.  This is currently under investigation.

\appendix

\section{Some results in Number Theory}
\label{Appendix}
\vspace{0.1in}

In this section, we briefly mention some fundamental results of Diophantine's approximation. The theorem of Khintchine (Theorem \ref{manyak}) plays an important role to determine the Lebesque measure of sets investigated in this paper.

Let $f:\mathbb{N}\to \mathbb{R}^+$ be called an approximation function   if $$\mathop {\lim }\limits_{\tilde q \to \infty } f(\tilde q)=0.$$ A real number  $\zeta$ is $f-approximable$ if $\zeta$ satisfies
\begin{eqnarray}\left|~\zeta-\frac{\tilde p}{\tilde q}~\right|<f(\tilde q)\label{app}
\end{eqnarray}
 for infinitely many rational numbers $\frac{\tilde p}{\tilde q}.$ Let $P(f)$ be the set of all $f-$approximable numbers. We recall the following theorem to find the measure of sets of type $P(f).$
\vspace{0.1in}

\begin{theorem} [Khintchine's theorem] \label{manyak}\cite[Page 4]{Bernik} \label{Khintchine} Let $\mu$ be the Lebesque measure. Then
\begin{eqnarray} \mu (P(f))=\left\{
                   \begin{array}{ll}
                     0, & {\text{if }} \quad \sum\limits_{\tilde q\in \mathbb{N}} {\tilde q f(\tilde q)}<\infty, \\
                     {\text{full}}, & {\text if} \quad \tilde q f(\tilde q) {\text{~is ~nonincreasing ~ and}} ~ \sum\limits_{\tilde q\in \mathbb{N}} {\tilde q f(\tilde q)}=\infty.
                   \end{array}
                 \right.
\end{eqnarray}
\end{theorem}
\vspace{0.1in}

 Dirichlet's theorem \cite{Cassal} states that every irrational number can be approximated to the order 2. The following theorem from \cite{Scott} is a special case of Dirichlet's theorem:
\vspace{0.1in}

\begin{theorem}\label{sac} Let $\zeta\in \mathbb{R}-\mathbb{Q}.$ Then there exists a constant $C\ge 1,$ and increasing sequences of coprime odd integers $\{\tilde p_j\},\{\tilde q_j\}$   satisfying the asymptotic relation
\begin{eqnarray}\label{irr-1}\left|~\zeta - \frac{\tilde p_j}{\tilde q_j}~\right|\le \frac{C}{{\tilde q_j}^{2}},\quad j\to\infty.\end{eqnarray}
\end{theorem}
\vspace{0.1in}

It obvious by Theorem \ref{Khintchine} that the set $ \mathbb{R}-\mathbb{Q}$ is uncountable and it has a full Lebesque measure.
\vspace{0.1in}
\begin{definition}
A real number $\zeta$ is a Liouville's number if for every $m\in \mathbb{N}$ there exists $\frac{\tilde p_m}{\tilde q_m}$ with $p_m, q_m \in\mathbb{Z}$ such that  $$\left|~\zeta-\frac{\tilde p_m}{\tilde q_m}~\right|<\frac{1}{\tilde q_m^m}.$$
\end{definition}
\vspace{0.1in}

It is proved that any Liouville's number is transcendental. Theorem \ref{Khintchine} implies that the set of Liouville's numbers is of Lebesque measure zero.

\vspace{0.1in}
\begin{definition} A real number $\zeta$ is an algebraic number if it is a root of a polynomial equation
$$a_n x^n+ a_{n-1}x^{n-1}+ \ldots + a_1 x + a_0=0$$
with each $a_i \in \mathbb{Z},$ and at least one of $a_i$ is non-zero. A number which is not algebraic is called transcendental.
\end{definition}

\vspace{0.1in}
 Now we give the following results of Diophantine's approximations:
 \vspace{0.1in}

\begin{theorem}\label{irr2} There exists a set $\tilde {\mathbb{Q}}$ such that if $\zeta\in \mathbb{R}-\tilde {\mathbb{Q}},$ then for every $\varepsilon>0$ there are infinitely many  $\frac{\tilde p}{\tilde q}\in\mathbb{Q}$ and a constant $C_\zeta>0$ such that
\begin{eqnarray}\left|~\zeta-\frac{\tilde p}{\tilde q}~\right|\ge \frac{C_\zeta}{\tilde q^{2+\varepsilon}}.\label{aptal}\end{eqnarray}
Moroever, $\mu (\tilde {\mathbb{Q}})=0.$
\end{theorem}
\vspace{0.1in}

    \textbf{Proof:}
We  know that the irrational algebraic numbers belong to $\tilde {\mathbb{Q}}$  by Roth's theorem  (Page 103, \cite{Cassal}). Therefore $\tilde {\mathbb{Q}}$ is not empty.  We proceed to the second part of the lemma. The first part of the theorem implies that  if  $\zeta\in \tilde {\mathbb Q}$ then for all $C_\zeta>0,$ the inequality $\left|~\zeta-\frac{\tilde p}{\tilde q}~\right|< \frac{C_\zeta}{\tilde q^{2+\varepsilon}}$  holds for some $\frac{\tilde p}{\tilde q}\in\mathbb{Q}.$ Now define the set
$$\tilde {\mathbb Q}_{\varepsilon}=\left\{\zeta\in \mathbb{R}~:~  \left|~~\zeta-\frac{\tilde p}{\tilde q}~\right|< \frac{C_\zeta}{\tilde q^{2+\varepsilon}}  ~ {\text{~ for ~infinitely ~many}} ~\frac{\tilde p}{\tilde q}\in\mathbb{Q} \right\}.$$
By the  notation of Theorem \ref{Khintchine}, choose $f(\tilde q)=\frac{C_\zeta}{\tilde q^{2+\varepsilon}}$ so that $\tilde q f(\tilde q)$ is nonincreasing and $\sum\limits_{\tilde q\in \mathbb{N}}\frac{C_\zeta}{\tilde q^{1+\varepsilon}}<\infty.$  By Theorem \ref{Khintchine}, $\mu (\tilde {\mathbb Q}_{\varepsilon})=0. $   Now we prove  $\tilde {\mathbb{Q}}\subset \tilde {\mathbb Q}_{\varepsilon}$ by contradiction. Assume that $\zeta\notin  \tilde {\mathbb Q}_{\varepsilon},$ i.e. there are finitely many  rationals $\left\{\frac{p_i}{q_i}\right\}_{i=1,\cdots, N}$ such that
 $$\left|~\zeta-\frac{p_i}{q_i}~\right|< \frac{C_\zeta}{q_i^{2+\varepsilon}} ~~for~~ i=1,\cdots, N, ~~{and }~~\left|~\zeta-\frac{\tilde p}{\tilde q}\right|\ge  \frac{C_\zeta}{q^{2+\varepsilon}} ~~ {\text{for}} ~~ \frac{\tilde p}{\tilde q}\notin
\bigcup\limits_{i = 1}^N {\left\{\frac{\tilde p_i}{\tilde q_i}\right\}}.$$
 The last inequality implies that $\zeta\in\mathbb{R}-\mathbb{Q}.$ 
 This implies that the set $\mathbb{R}-\tilde {\mathbb{Q}}$ has a full Lebesque measure. $\square$

\vspace{0.1in}

Now define the set $\tilde{\tilde {\mathbb{Q}}}$ by
\begin{eqnarray}\tilde{\tilde {\mathbb{Q}}}=\left\{\zeta\in \mathbb{R}~:~ \exists C>0, ~~ \left|~\zeta-\frac{\tilde p}{\tilde q}~\right|\ge  \frac{C}{\tilde q^{2}}  {\text{~ for ~infinitely ~many}} ~\frac{\tilde p}{\tilde q}\in\mathbb{Q} \right\}\label{yeni}.\end{eqnarray}
 If we consider numbers $\zeta\in\mathbb{R}$ whose the partial quotients satisfy $|a_k|<C(\zeta)$ for all $k\in\mathbb{N}$ in its continued fraction expansion
$$\zeta=[a_0; a_1, a_2, . . . ]=a_0+\frac{1}{a_1+\frac{1}{a_2+  \ddots}},$$
then $\zeta\in \tilde{\tilde {\mathbb{Q}}}.$ By Liouville's theorem (Page 128, \cite{Miller}), $\tilde{\tilde {\mathbb{Q}}}$ also contains all quadratic irrational numbers (the roots of an algebraic polynomial of
degree $2$). Therefore the set is uncountable.
\vspace{0.1in}

\begin{lemma}
The set $\tilde{\tilde {\mathbb{Q}}}$ has a Lebesgue measure zero.
\end{lemma}
\vspace{0.1in}

\textbf{Proof:} Define the set $F_m$ by
$$F_m=\left\{\zeta\in \mathbb{R}~:~  \left|~\zeta-\frac{\tilde p}{\tilde q}~\right|< \frac{C}{m\tilde  q^{2}}  ~ {\text{~ for ~infinitely ~many}}~\frac{\tilde p}{\tilde q}\in\mathbb{Q} \right\}.$$
Then $F_m$ has a full Lebesque measure by Theorem \ref{Khintchine}, i.e. $f(\tilde q)=\frac{C}{m\tilde q^2},$ and $\sum\limits_{\tilde q\in \mathbb{N}} {\tilde q f(\tilde q)}=\infty.$
Now consider the set  $\bigcap\limits_{m\in\mathbb{N}} F_m.$ This set is the countable intersection of sets $F_m,$ and each $F_m$ has full Lebesgue measure. Therefore $\mu \left(\bigcap\limits_{m\in\mathbb{N}} F_m\right)$ has full Lebesque measure.  Since  $\tilde{\tilde {\mathbb{Q}}}=\mathbb{R}-\bigcap\limits_{m\in\mathbb{N}} F_m,$ then $\mu (\tilde{\tilde {\mathbb{Q}}})=0.~\square$

\begin{acknowledgements}
I would like to thank to Prof. Kirsten Morris and Prof. Sergei Avdonin for the fruitful discussions and suggestions to finalize this paper.\end{acknowledgements}


\begin{thebibliography}{}

\bibitem{Alabou1} F. Alabau; P. Cannarsa, V. Komornik (2002){Indirect internal stabilization of weakly coupled evolution equations,} J. Evol. Equ. 2- 2, pp. 127--150.

\bibitem{Alabou2}  F. Alabau-Boussouira, M. Léautaud (2012) {Indirect stabilization of locally coupled wave-type systems,} ESAIM COCV 18-2, pp. 548--582.

\bibitem{Alabou3} F. Alabau-Boussouira (2003) { A two-level energy method for indirect boundary observability and controllability of weakly coupled hyperbolic systems} SIAM J. Control Optim. 42-3, pp. 871--906.

\bibitem{T1} F. Ammar-Khodja, A. Benabdallah, M. González-Burgos, L. de Teresa (2011) {Recent results on the controllability of linear coupled parabolic problems: a survey,} Math. Control Relat. Fields 1--3, pp. 267--306.

\bibitem{T2} F. Ammar Khodja, A. Benabdallah, M. González-Burgos, L. de Teresa (2013) {A new relation between the condensation index of complex sequences and the null controllability of parabolic systems,} C. R. Math. Acad. Sci. Paris 351/19--20, pp.  743--746.


\bibitem{Ammari-T} K. Ammari, M. Tucsnak (2001) {{Stabilization of second order evolution equations by a  class of unbounded feedbacks, }} ESAIM COCV 6, pp. 361--386.


\bibitem{AvdoninII}    S. Avdonin, W. Moran (2001) {{  Ingham-type inequalities and Riesz bases of divided differences,}}  Int. J. Appl. Math. Comput. Sci. 11--4, pp. 803--820.

\bibitem{AvdoninIII} S. Avdonin, W. Moran (2002) {{Riesz bases of exponentials and divided differences,}} St. Petersburg Math. J. 13--3, pp. 339--351.

\bibitem{Komornik-P2} C. Baiocchi, V. Komornik,  P. Loreti (2002) { Ingham-Beurling type theorems with weakened gap conditions,} Acta Math. Hungar., 97, pp. 55--95.


\bibitem{Banks-Smith} H.T. Banks, R.C. Smith, Y. Wang (1996) { Smart material structures: Modelling, Estimation and Control}, Mason, Paris, 1996.

\bibitem{Baz} Baz, A. (1996) { Active Constrained Layer Damping,} U.S. Patent \# 5,485,053.

\bibitem{Begout-Soria} P. B\'{e}gout, F. Soria (2007) { A generalized interpolation inequality and its application to the stabilization of damped equations,} J, Differ. Equations 240-2, pp. 324--356

\bibitem{Bernik}  V.I. Bernik, M.M. Dodson (1999) {{Metric Diophantine Approximation on Manifolds,}} Cambridge University Press,  Cambridge.

\bibitem{Cassal} J.W. Cassels (1966) { {An Introduction to Diophantine Approximation,}} Cambridge University Press, Cambridge.

\bibitem{C-Z} C. Castro, E. Zuazua (1998) {Boundary controllability of a hybrid system consisting in two
flexible beams
connected by a point mass,} SIAM J. Control Optim. 36-5, pp. 1576--1595.

\bibitem{C-W} R. Curtain, G. Weiss (1989) {{ Well posedness of triples of operators (in the sense of linear systems theory)}}, in: F. Kappel, K. Kunisch, W. Schappacher (Eds.), Control and Estimation of Distributed Parameter Systems, BirkhVauser, Basel 91, pp. 41–-59.


\bibitem{Zuazua}  R. Dager, E. Zuazua,{{Wave propagation, Observation and Control in 1--d Flexible Multi-structures,}} Springer,  2006.

\bibitem{Dest} Ph. Destuynder, I. Legrain, L. Castel, N. Richard (1992) { Theoretical, numerical and experimental discussion of the use of piezoelectric devices for control-structure interaction,} European J. Mech. A Solids, 11, pp. 181–-213.


\bibitem{Hansen3}
{  S.W. Hansen} (2004)
{ Several Related Models for Multilayer Sandwich Plates,}
{  Mathematical Models \& Methods in Applied
Sciences}  (14-8), pp. 1103-1132.

\bibitem{O-Hansen1}{S.W. Hansen, A.\"{O}. \"{O}zer} (2010) \newblock{ Exact boundary controllability of an abstract Mead-Marcus Sandwich beam model,} \newblock{ the Proceedings of 53rd IEEE Conf. on Decision \& Control}, Atlanta, USA, pp. 2578-2583.


\bibitem{J-T}  S. Jaffard, M. Tucsnak (1997) { Regularity of plate equations with control concentrated in interior curves,} Proc. Roy. Soc. Edinburg Sect. A 127, pp. 1005--1025.


\bibitem{Jaffard} S. Jaffard, M. Tucsnak, E. Zuazua (1998) { Singular internal stabilization of the wave equation,} J. Differential Equations 145-1, pp. 184--215.
\bibitem{K-M-M} B. Kapitonov, B. Miara, and G.P. Menzala (2007) {{Boundary Observation and Exact Control of a Quasi-electrostatic Piezoelectric System in Multilayered Media,}} SIAM J. Control Optim. 46-3, pp. 1080--1097.

\bibitem{Komornik-P} V. Komornik, P. Loreti (2005) { Fourier Series in Control Theory }, Springer-Verlag, New York.


\bibitem{Lions}
J.L. Lions (1988) {{Exact Controllability, stabilization and perturbations for distributed parameter systems.}}
SIAM Rev.  30-1, pp. 1--68.

\bibitem{T3} F. Luca, L. de Teresa (2013) {Control of coupled parabolic systems and Diophantine approximations},  SeMA J. 61, pp. 1--17.


   \bibitem{Miller} S.J. Miller and R. Takloo-Bighash (2006) {{ An Invitation to Modern Number Theory,}} Princeton University Press,
Princeton, NJ, 2006.

\bibitem{O-M}  K.A. Morris, A.\"{O}. \"{O}zer, (2014) {{Comparison  of stabilization of  current-actuated  and voltage-actuated piezoelectric beams,}} the Proceedings of 53rd IEEE Conf. on Decision \& Control, Los Angeles, USA, pp. 571--576.


\bibitem{cdcpaper} K.A. Morris, A.\"{O}. \"{O}zer (2013) { Strong stabilization of  piezoelectric beams with magnetic effects}, the Proceedings of 52nd IEEE Conf. on Decision \& Control, Firenze, Italy, pp. 3014--3019.

\bibitem{O-M1}  K.A. Morris, A.\"{O}. \"{O}zer (2014) {{Modeling and stabilizability of voltage-actuated piezoelectric beams with magnetic effects,}}  SIAM J. Control Optim. (52--4), pp. 2371--2398.

\bibitem{O-Hansen3} {A.\"{O}. \"{O}zer, S.W. Hansen} (2014)  \newblock{ Exact boundary controllability results for a multilayer Rao-Nakra sandwich beam}, SIAM J. Control Optim. (52-2), pp. 1314--1337.


\bibitem{O-Hansen4} {A.\"{O}. \"{O}zer, S.W. Hansen} (2013) \newblock{Uniform stabilization of a multi-layer Rao-Nakra sandwich beam},  Evolution Equations and Control Theory, (2-4), pp. 195--210.



\bibitem {accpaper} A.\"{O}. \"{O}zer, K. A. Morris (2014) { Modeling an elastic beam with piezoelectric patches by including magnetic effects,} the Proceedings of the American Control Conference, Portland, USA, pp. 1045-1050.



\bibitem{Rogacheva} N. Rogacheva (1994) {{ The Theory of Piezoelectric Shells and Plates}}, Boca Raton, FL: CRC Press.


\bibitem{Ronkanen} P. Ronkanen, P. Kallio, M. Vilkko, H.N. Koivo (2011)  { Displacement Control of Piezoelectric Actuators Using Current and Voltage,} IEEE/ASME Trans. Mechatronics 16-1, pp. 160--166.


\bibitem{Russell} D.L. Russell (1986) {{ The Dirichlet–Neumann boundary control problem associated with Maxwell's
equations in a cylindrical region,}} {{SIAM J. Control Optim.,}} 24, pp. 199–-229.


\bibitem{Scott} W.T. Scott (1940) {{Approximation to real irrationals by certain classes of rational fractions,}} Bull. Amer. Math. Soc. 46, pp. 124--129.

\bibitem{Smith} R.C. Smith (2005) { Smart Material Systems}, Society for
Industrial and Applied Mathematics.

\bibitem{Tiersten} H.F. Tiersten (1969) { {Linear Piezoelectric Plate Vibrations} }, Plenum Press, New York.
\bibitem{Triebel} H. Triebel (1978), Interpolation Theory, Function Spaces, Differential Operators.
North-Holland, Amsterdam.
\bibitem{Tucsnak-a} M. Tucsnak (1996) { Regularity and exact controllability for a beam with piezoelectric actuator, } SIAM J. Cont. Optim. 34, pp. 922--930.
\bibitem{Tucsnak-b} M. Tucsnak (1996) { Control of plate vibrations by means of piezoelectric actuators, } Discrete Contin. Dynam. Systems 2, pp.
 281--293.
\bibitem{Weiss-Tucsnak1} M. Tucsnak, G. Weiss (2001) {{Simultaneous controllability in sharp time for two elastic strings}}, ESAIM: Cont. Optim. Calc. Var. 6, pp. 259--273.
\bibitem{Weiss-Tucsnak} M. Tucsnak, G. Weiss (2009) {{Observation and Control for Operator Semigroups}}, Birkhäuser Verlag, Basel.

\bibitem{Weiss} M. Tucsnak, G. Weiss (2003) {{How to get a conservative well-posed linear system out of thin air, Part II: controllability and stability,}} SIAM J. Cont. Optim. 42-3, pp. 907-935.



\bibitem{Tzou} H.S. Tzou (1993) {{Piezoelectric shells, Solid Mechanics and Its applications 19,}} Kluwer Academic, The Netherlands.

\bibitem{Ulrich} D. Ulrich (1980) {Divided Differences and Systems of Nonharmonic Fourier Series,} Proc. of the
Amer. Math. Soc., 80--1, pp. 47–-57.

\bibitem{Weiss-a} G. Weiss, M. Tucsnak (2003) {{How to get a conservative well-posed linear system out of thin air, Part I. Well-posedness and energy balance,}} ESAIM:  Control, Optimization, Calculus of Variations 9, pp. 247--273.

    \bibitem{Yang} J. Yang (2005) {{An Introduction to the Theory of Piezoelectricity,}} Springer, New York.

    \bibitem{Yang1} J. Yang (2006) {{A review of a few topics in piezoelectricity,}} Appl. Mech. Rev. 59, pp. 335–-345.
%
%
\end{thebibliography}


\end{document}